\documentclass[aap,MSNbibl,citesort,seceqn,dvips]{arximspdf}
\usepackage[vtexbibl]{}

\doi{10.1214/09-AAP625}
\volume{20}
\issue{2}
\pubyear{2010}
\firstpage{495}
\lastpage{521}

\makeatletter

\newtheorem{theorem}{Theorem}[section]
\newtheorem{lem}[theorem]{Lemma}
\newtheorem{Claim}{Claim}
\newproclaim{Case}{Case}
\newproclaim{Fact}{Fact}
\newtheorem{prop}[theorem]{Proposition}
\newproclaim{definition}[theorem]{Definition}
\newproclaim{example}[theorem]{Example}
\newproclaim{rmk}[theorem]{Remark}
\newproclaim{Acknowledgments}{Acknowledgments}
\makeatother

\begin{document}
\begin{frontmatter}

\title{A Birthday Paradox for Markov chains with an optimal bound for
collision in the Pollard Rho~algorithm for discrete logarithm}
\runtitle{A Birthday Paradox for Pollard Rho}

\begin{aug}
\author[a]{\fnms{Jeong Han} \snm{Kim}\ead[label=e1]{jehkim@yonsei.ac.kr}\ead[label=e5]{jehkim@nims.re.kr}\thanksref{t1}},
\author[b]{\fnms{Ravi} \snm{Montenegro}\corref{}\ead[label=e2]{ravi\_montenegro@uml.edu}},\\
\author[c]{\fnms{Yuval} \snm{Peres}\ead[label=e3]{peres@microsoft.com}\thanksref{t3}}
\and
\author[d]{\fnms{Prasad} \snm{Tetali}\ead[label=e4]{tetali@math.gatech.edu}\thanksref{t4}}
\runauthor{Kim, Montenegro, Peres and Tetali}
\thankstext{t1}{Supported by  Korea Science and Engineering
Foundation (KOSEF) Grant funded by the Korean government (MOST)
R16-2007-075-01000-0.}
\thankstext{t3}{Supported in part by NSF Grant DMS-06-05166.}
\thankstext{t4}{Supported in part by NSF Grants DMS-04-01239, DMS-07-01043.}
\affiliation{Yonsei University, University of Massachusetts Lowell,
Microsoft Research and Georgia Institute of Technology}

\address[a]{J. H. Kim\\ Department of Mathematics\\ Yonsei
University\\ Seoul, 120-749\\ Korea\\ and\\ National Institute for\\
\quad Mathematical Sciences\\ 628 Daeduk Boulevard\\ Daejeon, 305-340\\
Korea\\
\printead{e1}\\
\phantom{E-mail: }\printead*{e5}}
\address[b]{R. Montenegro\\ Department of Mathematical Sciences\\
University of Massachusetts Lowell\\ Lowell, Massachusetts 01854\\ USA\\
\printead{e2}}
\address[c]{Y. Peres\\ Microsoft Research\\ One Microsoft Way\\ Redmond, Washington 98052\\ USA\\
\printead{e3}}
\address[d]{P. Tetali\\ School of Mathematics\\ Georgia Institute of
Technology\\ Atlanta, Georgia 30332\\ USA\\
\printead{e4}}
\end{aug}

% HISTORY:
\received{\smonth{4} \syear{2008}}
\revised{\smonth{6} \syear{2009}}

% ABSTRACT
%
\begin{abstract}
We show a Birthday Paradox for self-intersections of Markov chains with
uniform stationary distribution.
As an application, we analyze Pollard's Rho algorithm for finding the
discrete logarithm in a cyclic group
$G$ and find that if the partition in the algorithm is given by a
random oracle, then with high probability a
collision occurs in $\Theta(\sqrt{|G|})$ steps. Moreover, for the
parallelized distinguished points algorithm
on $J$ processors we find that $\Theta(\sqrt{|G|}/J)$ steps suffices.
These are the first proofs of the correct
order bounds which do not assume that every step of the algorithm
produces an i.i.d. sample from $G$.
\end{abstract}

% KEYWORDS
%


\end{frontmatter}

%s1 ###
\section{Introduction}

The Birthday Paradox states that if $C\sqrt{N}$ items are sampled
uniformly at random with replacement from a set of $N$ items, then, for
large~$C$ with high probability, some items will be chosen twice. This
can be interpreted as a statement that with high probability, a Markov
chain on the complete graph $K_N$ with transitions $P(i,j)=1/N$ will
intersect its past in $C\sqrt{N}$ steps; we refer to such a
self-intersection as a \textit{collision} and say the ``\textit{collision
time}'' is $O(\sqrt{N})$.
Miller and Venkatesan generalized this in \cite{MV061} by showing
that for a general Markov chain, the collision time is bounded
by $O(\sqrt{N} T_s(1/2))$ where $T_s(\varepsilon)=\min\{n\dvtx \forall
u,v\in V, P^n(u,v) \ge(1-\varepsilon)\pi(v)\}$ measures the time
required for the $n$-step distribution to assign every state a suitable
multiple of its stationary probability. Kim, Montenegro and Tetali
\cite{KMT071} further improved the bound on collision time to
$O(\sqrt{N T_s(1/2)})$. In contrast, while this shows the average
path to be quickly self-intersecting, Pak \cite{Pak021} has shown
that undirected regular graphs of large degree have a nonintersecting
path of length $N/(32T_s(1/2))$.
%The walk on $K_N$ has $T_s(1/2)=O(1)$, and so these extend the
%Birthday Paradox.

The motivation of \cite{MV061,KMT071} was to study the collision
time for a Markov chain involved in Pollard's Rho algorithm for finding
the discrete logarithm on a cyclic group $G$ of prime order $N=|G|\neq
2$. For this walk, $T_s(1/2)=\Omega(\log N)$, and so the results of
\cite{MV061,KMT071} are insufficient to show the widely believed
$\Theta(\sqrt{N})$ collision time for this walk. In this paper we
improve upon these bounds and show that if a finite ergodic
Markov
chain has uniform stationary\vspace{1pt} distribution over $N$ states, then
$O(\sqrt{N})$ steps suffice for a collision to occur as long as the
relative-pointwise distance ($L_\infty$ of the densities of the
current and the stationary distribution) drops steadily \textit{early} in
the random walk; it turns out that the precise mixing time is largely,
although not entirely, unimportant. See Theorem \ref{theorem:birthday} for
a precise statement.
This is then applied to the Rho walk to give the first proof of
collision in $\Theta(\sqrt{N})$ steps, matching Shoup's lower bound \cite{Shoup97}
on time required for any probabilistic generic algorithm to solve this problem,
 and to van Oorschot and
Wiener's \cite{VW99} parallel version of the algorithm on $J$
processors to prove collision in $\Theta(\sqrt{N}/J)$ steps.

We note here that it is also well known (see, e.g., \cite{AF}, Section
4.1) that a random walk of length $L$
contains roughly $L\lambda$ samples from the
stationary measure (of the Markov chain) where $\lambda$ is the
spectral gap of the chain. This yields another estimate on collision
time for a Markov chain which is also of a multiplicative nature
(namely, $\sqrt{N}$ \textit{times}
a function of the mixing time) as in \cite{MV061,KMT071}.
A main point of the present work
is to establish sufficient criteria under which the collision time has
an \textit{additive}
bound: $C\sqrt{N}$ plus an estimate on the mixing time. While the Rho algorithm
provided the main motivation for the present work, we find the more
general Birthday Paradox result to be of independent interest, and, as
such, expect to have other applications in the future.

A bit of detail about the Pollard Rho algorithm is in order.
The classical discrete logarithm problem on a cyclic group deals with
computing the exponents, given the generator of the group; more
precisely, given a generator $g$ of a cyclic group $G$ and an element
$h=g^x$, one would like to compute $x$ efficiently.
Due to its presumed computational difficulty, the problem figures
prominently in various cryptosystems, including the Diffie--Hellman key
exchange, El Gamal system and
elliptic curve cryptosystems. About 30 years ago, J. M. Pollard suggested
algorithms to help solve both factoring large integers \cite{Pol75}
and the discrete logarithm problem \cite{Pol78}.
While the algorithms are of much interest in computational number
theory and cryptography, there has been little work on rigorous
analysis. We refer the reader to \cite{MV061} and other existing
literature (e.g., \cite{Tes01,CP05}) for further cryptographic and
number-theoretical motivation for the discrete logarithm
problem.

A standard variant of the classical Pollard Rho algorithm for finding
discrete logarithms can be described using a Markov chain on a cyclic
group $G$.
While there has been no rigorous proof of rapid mixing of this Markov
chain of order $O(\log^c |G|)$ until recently, Miller and Venkatesan
\cite{MV061} gave a proof of mixing of order $O(\log^3|G|)$ steps
and collision time of $O(\sqrt{|G|}\log^3|G|)$, and Kim, Montenegro and Tetali
\cite{KMT071} showed mixing of order $O(\log|G| \log\log|G|)$ and
collision time of $O(\sqrt{|G|\log|G| \log\log|G|})$. In this
paper we give the first proof of the correct $\Theta(\sqrt{|G|})$
collision time. By recent results of Miller and Venkatesan \cite{MV081}
this collision will be nondegenerate and will solve the discrete
logarithm problem with probability $1-o(1)$ for almost every prime
order $|G|$, if the start point of the algorithm is chosen at random or
if there is no collision in the first $O(\log|G| \log\log|G|)$ steps.

The paper proceeds as follows. Section \ref{sec:prelim} contains some
preliminaries, primarily an introduction to the Pollard Rho algorithm
and a simple multiplicative bound
on the collision time in terms of the mixing time. The more general
Birthday Paradox for Markov chains with uniform stationary distribution
is shown in Section \ref{sec:collision}. In Section \ref{sec:mixing}
we bound the appropriate constants for the Rho walk and show the
optimal collision time. We finish in Section \ref{sec:parallel} by
proving similar results for the distinguished points method of
parallelizing the algorithm.

%*********************** Preliminaries *******************************

%s2 ###
\section{Preliminaries} \label{sec:prelim}

Our intent in generalizing the Birthday Paradox was to bound the
collision time of the Pollard Rho algorithm for discrete logarithm. As
such, we briefly introduce the algorithm here. Throughout the analysis
in the following sections, we assume that the size $N=|G|$ of the
cyclic group on which the random walk is performed is odd. Indeed there
is a standard reduction (see \cite{Pom} for a very readable account
and \cite{PohHel78} for a classical reference) justifying the fact
that it suffices to study the discrete logarithm problem on cyclic
groups of \textit{prime} order.

Suppose $g$ is a generator of $G$, that is, $G=\{g^i\}_{i=0}^{N-1}$.
Given $h\in G$, the discrete logarithm problem asks us
to find $x$ such that $g^x=h$. Pollard suggested an algorithm on
$\mathbb{Z}
_N^{\times}$ based on a random walk and the Birthday
Paradox. A common extension of his idea to groups of prime order is to
start with a partition of $G$ into sets $S_1$, $S_2$,
$S_3$ of roughly equal sizes, and define an iterating function $F\dvtx
G\to G$ by $F(y)=gy$ if $y\in S_1$, $F(y)=hy=g^x y$ if $y\in S_2$ and
$F(y)=y^2$ if $y\in S_3$.
Then consider the walk $y_{i+1}=F(y_i)$. If this walk passes through
the same
state twice, say $g^{a+xb}=g^{\alpha+x\beta}$, then $g^{a-\alpha
}=g^{x(\beta-b)}$, and so $a-\alpha\equiv x(\beta-b)\operatorname{mod} N$
and $x\equiv(a-\alpha)(\beta-b)^{-1} \operatorname{mod} N$ which determines $x$ as
long as $(\beta-b,N)=1$ (the nondegenerate case). Hence, if we define a
\textit{collision} to be the event that the walk passes over the same
group element twice, then the first time there is a
collision it might be possible to determine the discrete logarithm.

To estimate the running time until a collision, one heuristic is to
treat $F$ as if it outputs uniformly random group elements. By the
Birthday Paradox, if $O(\sqrt{|G|})$ group elements are chosen
uniformly at random, then there is a high probability that two of these
are the same. Teske \cite{Tes98} has given experimental evidence that
the time until a collision is slower than what would be expected by an
independent uniform random process.
We analyze instead the actual
Markov chain in which it is assumed only that each $y\in G$ is assigned
independently and at random to a partition $S_1$, $S_2$ or $S_3$. In
this case, although the iterating function $F$ described earlier is
deterministic, because the
partition of $G$ was randomly chosen, then the walk is equivalent to a
Markov chain (i.e., a random walk), at least until
the walk visits a previously visited state and a collision occurs. The
problem is then one of considering a walk on the exponent of $g$, that
is, a walk $P$ on the cycle $\mathbb{Z}_N$ with transitions $P(u,u+1)=
P(u,u+x)= P(u,2u)= 1/3$.

\begin{rmk}
By assuming each $y\in G$ is assigned independently and at random to a
partition, we have eliminated one of the key features of the Pollard
Rho algorithm, space efficiency. However, if the partitions are given
by a hash function, $f\dvtx (G,N)\to\{1,2,3\}$, which is sufficiently
pseudo-random, then we might expect behavior similar to the model with
random partitions.
\end{rmk}

\begin{rmk}
While we are studying the time until a collision occurs, there is no
guarantee that the first collision will be nondegenerate. If the first
collision is degenerate then so will be all collisions as the algorithm
becomes deterministic after the first collision.
\end{rmk}

As mentioned in the introduction, we first recall a simple
multiplicative bound on
collision time from \cite{KMT071}. The following proposition relates
$T_s(1/2)$ to the time until a collision occurs for any Markov chain
$P$ with uniform distribution on $G$ as the stationary distribution.
%Recall that the event of revisiting an already visited state is called
%a collision.

\begin{prop} \label{prop:oldresult}
With the above definitions, a collision occurs after
\[
T_s(1/2) + 2 \sqrt{2c |G| T_s(1/2)}
\]
steps with probability at least $1-e^{-c}$, for any $c>0$.
\end{prop}

\begin{pf}
Let $S$ denote the first $\lceil\sqrt{2c |G| T_s(1/2)}
\rceil$ states visited by the walk. If two of these states are the
same then a collision has occurred, so assume all states are distinct.
Even if we only check for collisions every $T_s(1/2)$ steps, the chance
that no collision occurs in the next $tT_s(1/2)$ steps (so consider $t$
semi-random states) is then at most
\[
\biggl(1-\frac12\frac{|S|}{|G|}\biggr)^t
\leq\biggl(1-\sqrt{\frac{c T_s(1/2)}{2|G|}}\biggr)^t
\leq\exp\biggl(-t \sqrt{\frac{c T_s(1/2)}{2|G|}}\biggr) .
\]
When $t=\lceil\sqrt{\frac{2c |G|}{T_s(1/2)}}\rceil$,
this is at most $e^{-c}$, as desired, and so at most
\[
\bigl\lceil\sqrt{2c |G| T_s(1/2)}\bigr\rceil-1
+ \biggl\lceil\sqrt{\frac{2c |G|}{T_s(1/2)}}\biggr\rceil T_s(1/2)
\]
steps are required for a collision to occur with probability at least
$1-e^{-c}$.
\end{pf}

Obtaining a more refined additive bound on collision time will be the
focus of the next
section. While the proof can be seen as another application of the
well-known second moment method, it turns out that bounding the second
moment of the number of
collisions \textit{before} the mixing time is somewhat subtle. To handle
this, we use
an idea from \cite{LPS03}, who in turn credit their approach to \cite{GR91}.
%%should we say more here about what LPS prove regarding the variance of
%%the number of collisions of two independent copies of a vertex
%transitive chain?!

%*********************** Collision Times *****************************

%s3 ###
\section{Collision time} \label{sec:collision}

Consider a finite ergodic Markov chain $P$ with uniform stationary
distribution $U$ (i.e., doubly stochastic), state space $\Omega$ of
cardinality $N=|\Omega|$, and let $X_0,X_1,\ldots$ denote a
particular instance of the walk.
In this section we determine the number of steps of the walk required
to have a high probability that a ``collision'' has occurred, that is,
a self-intersection $X_i=X_j$ for some $i\neq j$.

A key notion when studying Markov chains is the mixing time, or the
time required until the probability of being at each state is suitably
close to its stationary probability.
\begin{definition}
The \textit{mixing time} $\tau(\varepsilon)$ of a Markov chain $P$ with
stationary distribution $U$ is given by
\[
\tau(\varepsilon)
= \min\{T\dvtx \forall u,v\in\Omega, (1-\varepsilon)U(v) \leq
P^T(u,v) \leq(1+\varepsilon)U(v)\} .
\]
\end{definition}

Now some notation.

Fix some $T\geq0$ and integer $\beta>0$. Let the indicator function
$\mathbf{1}_{\{X_i=X_j\}}$ equal one if $X_i=X_j$, and zero otherwise. Define
\[
S= \sum_{i=0}^{\beta\sqrt{N}} \sum_{j=i+2T}^{\beta\sqrt{N}+2T}
\mathbf{1}_{\{X_i=X_j\}}
\]
to be the number of times the walk intersects itself in $\beta\sqrt
{N}+2T$ steps where $i$ and $j$ are at least $2T$ steps apart. Also,
for $u,v \in\Omega$, let
\[
G_T(u,v) = \sum_{i=0}^T P^i(u,v)
\]
be the expected number of times a walk beginning at $u$ hits state $v$
in $T$ steps. Finally, let
\[
A_T = \max_u \sum_v G^2_T(u,v)\quad\mbox{and}\quad
A^*_T = \max_u \sum_v G^2_T(v,u) .
\]

To see the connection between these and the collision time, observe that
\begin{eqnarray*}
\sum_v G_T^2(u,v)
&=& \sum_v \Biggl(\sum_{i=0}^T\sum_{j=0}^TP^i(u,v)P^j(u,v) \Biggr)
\\
&=& \sum_{i=0}^T\sum_{j=0}^T \sum_v P^i(u,v)P^j(u,v)
\\
&=& \sum_{i=0}^T \sum_{j=0}^T \mathsf{Pr}(X_i=Y_j)
 \\
&=& \sum_{i=0}^T \sum_{j=0}^T E\bigl(\mathbf{1}_{\{X_i=Y_j\}}\bigr)
= E\sum_{i,j=0}^T \mathbf{1}_{\{X_i=Y_j\}},
\end{eqnarray*}
where $\{X_i\}, \{Y_j\}$ are i.i.d. copies of the chain, both having
started at $u$ at time 0, and $E$ denotes expectation. Hence $A_T$ is
the maximal expected number of collisions of two $T$-step i.i.d. walks
of $P$ starting at the same state $u$. Likewise, $A^*_T$ is the same
for the reversal $P^*$ where $P^*(u,v)=P(v,u)$ (recall that the
stationary distribution was assumed to be uniform).

The main result of this section is the following.

\begin{theorem}[(Birthday Paradox for Markov chains)] \label{theorem:birthday}
Consider a finite ergodic Markov chain with uniform stationary
distribution on a state space of size~$N$. Let $T$ be such that $\frac
{m}{N} \leq P^T(u,v) \leq\frac{M}{N}$ for some $m\leq1\leq M$ and
every pair of states $u,v$.
%The expected number of steps until a collision is at most
%$$
%4(\frac Mm)^2 (\sqrt{\frac{2N}{M}\max\{A_T,A_T^*\}}+T
%) ,
%$$
%and more specifically
After
\[
4c\biggl(\frac Mm\biggr)^2 \biggl(\sqrt{\frac{2N}{M}\max\{A_T,
A^*_T\}}+T\biggr)
\]
steps a collision occurs with a probability of at least $1-e^{-c}$, for
any $c\geq0$.
\end{theorem}

At the end of this section a slight strengthening of Theorem \ref
{theorem:birthday} is shown at the cost of a somewhat less intuitive bound.

In Example \ref{ex:pre-mixing}, near the end of this section, we
present an example to illustrate the need for the pre-mixing term $A_T$
in Theorem \ref{theorem:birthday}.
In contrast, very recently Nazarov and Peres \cite{NP08} proved a
general bound
for the birthday problem on any \textit{reversible} Markov chain on $N$
states: Suppose that the ratio of the stationary measures of any two
states is at most $A$. Then they
show that for any starting state, the expected time until the chain
visits a previously visited state is at most
$C\sqrt{N} + \log(A)$ for some universal constant $C$. In particular,
this implies an expected collision time of
$O(\sqrt{N})$ for the simple random walk on an undirected graph on $N$
vertices, and so the pre-mixing term is not necessary when considering
reversible walks.

Observe that if $A_T,A_T^*,m,M=\Theta(1)$ and $T=O(\sqrt{N})$, then
the collision time is $O(\sqrt{N})$ as in the standard Birthday
Paradox. By Lemma \ref{lem:quick-convergence}, for this to occur it
suffices that $P^T$ be sufficiently close to uniform after $T=o(\sqrt
{N})$ steps and that $P^j(u,v)=o(T^{-2})+d^j$ for all $u,v$, for $j\leq
T$ and some $d<1$. More generally, to upper bound $A_T$ and $A_T^*$ it
suffices to show that the maximum probability of being at a vertex
decreases quickly.

\begin{lem} \label{lem:quick-convergence}
If a finite ergodic Markov chain has uniform stationary distribution, then
\[
A_T,A_T^* \leq2\sum_{j=0}^T (j+1) \max_{u,v} P^j(u,v) .
\]
\end{lem}

\begin{pf}
If $u$ is such that equality occurs in the definition of $A_T$, then
\begin{eqnarray*}
A_T&=&\sum_v G^2_T(u,v) = \sum_{i=0}^T \sum_{j=0}^T \sum_v
P^i(u,v)P^j(u,v)
\\
&\leq& 2 \sum_{j=0}^T \sum_{i=0}^j \max_{y} P^j(u,y) \sum_v
P^i(u,v)
\\
&\leq& 2 \sum_{j=0}^T (j+1) \max_{y} P^j(u,y)  .
\end{eqnarray*}

The quantity $A_T^*$ plays the role of $A_T$ for the reversed chain,
and so the same bound holds for $A_T^*$ but with $\max
_{u,v}(P^*)^j(u,v)=\max_{u,v}P^j(v,u)=\max_{u,v}P^j(u,v)$.
\end{pf}

In particular, suppose $P^j(u,v)\leq c + d^j$ for every $u,v\in\Omega
$ and some $c,d\in[0,1)$. The sum
\begin{eqnarray*}
\sum_{j=0}^T (j+1)(c+d^j)
&=& c \frac{(T+1)(T+2)}{2} + \frac
{1-d^{T+1}-(T+1)d^{T+1}(1-d)}{(1-d)^2}
\\
&\leq& \bigl(1+o(1)\bigr)\frac{cT^2}{2} + \frac{1}{(1-d)^2} ,
\end{eqnarray*}
and so if $P^j(u,v) \leq o(T^{-2}) + d^j$ for every $u,v\in\Omega$,
then $A_T,A_T^*=\frac{2+o(1)}{(1-d)^2}$.

The proof of Theorem \ref{theorem:birthday} relies largely on the
following inequality which shows that the expected number of
self-intersections is large with low variance:

\begin{lem} \label{lem:expectations}
Under the conditions of Theorem \ref{theorem:birthday},
\begin{eqnarray*}
E[S] &\geq& \frac{m}{N} \pmatrix{\beta\sqrt{N}+2\cr 2},
\\
E[S^2] &\leq& \frac{M^2}{N^2}\pmatrix{\beta\sqrt{N}+2\cr 2}^2
\biggl(1+\frac{8\max\{A_T,A_T^*\}}{M\beta^2}\biggr).
\end{eqnarray*}
\end{lem}

\begin{pf*}{Proof of Theorem \ref{theorem:birthday}}
Recall the standard second moment bound: using Cauchy--Schwarz, we have that
\[
E [S] = E\bigl[ S \mathbf{1}_{\{S>0\}}\bigr]
\leq E[S^2 ]^{1/2} E\bigl[ \mathbf{1}_{\{S>0\}}\bigr]^{1/2},
\]
and hence
$\mathsf{Pr}[ S>0] \geq E [ S]^2/E [ S^2] .
$
If $\beta=2\sqrt{2\max\{A_T,A_T^*\}/M}$, then by Lemma \ref
{lem:expectations},
\[
\mathsf{Pr}[ S>0] \geq\frac{m^2/M^2}{1+{(8\max\{A_T,A_T^*\}
)/(M\beta^2)}}
\geq\frac{m^2}{2M^2}
\]
independent of the starting point. If no collision occurs in $\beta
\sqrt{N}+2T$ steps then $S=0$ as well, and so $\mathsf{Pr}[\mathit{no\
collision}]\leq
\mathsf{Pr}[S=0]\leq1-m^2/2M^2$. Hence, in $k(\beta\sqrt{N}+2T)$ steps,
%e1 ###
%
\begin{equation} \label{eqn:positive_prob}
\mathsf{Pr}[\mathit{no\ collision}]\leq(1-m^2/2M^2)^k\leq e^{-km^2/2M^2}.
\end{equation}
Taking $k=2cM^2/m^2$ completes the proof.
% of the final relation. For the expectation, let $s=\beta\sqrt{N}+2T$
%and $C$ be the number of steps until a collision. Then
%$$
%E[C] \leq\sum_{k=0}^{\infty} Pr[C>ks] s
% \leq\sum_{k=0}^{\infty} (1-\frac{m^2}{2M^2})^k s
% = 2(\frac Mm)^2 (\beta\sqrt{N}+2T) .
%$$
\end{pf*}

\begin{pf*}{Proof of Lemma \ref{lem:expectations}}
We will repeatedly use the relation that there are ${{\beta\sqrt
{N}+2}\choose {2}}$ choices for $i,j$ appearing in the summation for $S$, that
is, $0\leq i$ and $i+2T\leq j\leq\beta\sqrt{N}+2T$.

Now to the proof. The expectation $E[S]$ satisfies
%e2 ###
%
\begin{eqnarray}\label{eqn1}
E[S]
&=& E \sum_{i=0}^{\beta\sqrt{N}} \sum_{j=i+2T}^{\beta\sqrt{N}+2T}
\mathbf{1}_{\{X_i=X_j\}} \nonumber
\\[-8pt]\\[-8pt]
&=& \sum_{i=0}^{\beta\sqrt{N}} \sum_{j=i+2T}^{\beta\sqrt{N}+2T}
E\bigl[ \mathbf{1}_{\{X_i=X_j\}} \bigr]
\geq\pmatrix{\beta\sqrt{N}+2\cr 2}  \frac{m}{N}\nonumber
\end{eqnarray}
because if $j\geq i+T$, then
%e3 ###
%
\begin{eqnarray}\label{eqn2}
\quad \mathsf{Pr}(X_j=X_i)
= \sum_u \mathsf{Pr}(X_i=u)P^{j-i}(u,u)
\geq \sum_u \mathsf{Pr}(X_i=u)\frac{m}{N}
= \frac mN .
\end{eqnarray}
Similarly, $\mathsf{Pr}(X_j=X_i)\leq\frac{M}{N}$ when $j\geq i+T$.

Now for $E[S^2]$. Note that
\begin{eqnarray*}
E[S^2] &=& E\Biggl(\sum_{i=0}^{\beta\sqrt{N}} \sum_{j=i+2T}^{\beta
\sqrt{N}+2T} \mathbf{1}_{\{X_i=X_j\}}\Biggr)
\Biggl( \sum_{k=0}^{\beta\sqrt{N}} \sum_{l=k+2T}^{\beta\sqrt
{N}+2T} \mathbf{1}_{\{X_k=X_l\}}\Biggr)
\\
&=& \sum_{i=0}^{\beta\sqrt{N}} \sum_{k=0}^{\beta\sqrt{N}}
\sum_{j=i+2T}^{\beta\sqrt{N}+2T} \sum_{l=k+2T}^{\beta\sqrt{N}+2T}
\operatorname{Prob}(X_i=X_j, X_k=X_l)  .
\end{eqnarray*}
To evaluate this quadruple sum we break it into $3$ cases.

\begin{Case}\label{ca1}
Suppose $|j-l|\geq T$. Without loss, assume $l\geq j$ so that, in
particular, $l\geq\max\{i,j,k\}+T$. Then
%e4 ###
%
\begin{eqnarray}\label{eqn3}
&& \operatorname{Prob}(X_i=X_j,  X_k=X_l)  \nonumber
\\
&&\qquad = \operatorname{Prob}(X_i=X_j) \operatorname{Prob}(X_l=X_k \mid X_i=X_j)\nonumber
\\[-8pt]\\[-8pt]
&&\qquad \leq \operatorname{Prob}(X_i=X_j) \max_{u,v} \operatorname{Prob}\bigl(X_l=v \mid X_{\max\{i,j,k\}
}=u\bigr) \nonumber
\\
&&\qquad\leq \operatorname{Prob}(X_i=X_j) \frac{M}{N} \leq\biggl(\frac{M}{N}\biggr)^2.\nonumber
\end{eqnarray}
The first inequality holds because $\{X_t\}$ is a Markov chain and so
given $X_i,X_j,X_k$ the walk at any time $t\geq\max\{i,j,k\}$ depends
only on the state $X_{\max\{i,j,k\}}$.
\end{Case}

\begin{Case}\label{ca2}
Suppose $|i-k|\geq T$ and $|j-l|<T$. Without loss, assume $i\leq k$. If
$j\leq l$, then
%e5 ###
%
\begin{eqnarray}\label{eqn4}
&&\operatorname{Prob}(X_i=X_j, X_k=X_l)  \nonumber
\\
&&\qquad = \sum_{u,v} \operatorname{Prob}(X_i=u) P^{k-i}(u,v)P^{j-k}(v,u)P^{l-j}(u,v)
\\
&&\qquad \leq \sum_u \operatorname{Prob}(X_i=u) \frac{M}{N} \frac{M}{N} \sum_v
P^{l-j}(u,v) =\biggl(\frac{M}{N}\biggr)^2\nonumber
\end{eqnarray}
because $k\geq i+T$, $j\geq k+T$ and $\sum_v P^t(u,v)=1$ for any $t$
because $P$ and hence also $P^t$ are stochastic matrices. If, instead,
$l<j$ then essentially the same argument works, but with $\sum_v
P^t(v,u)=1$ because $P$ and hence also $P^t$ are doubly-stochastic.
\end{Case}

\begin{Case}\label{ca3}
Finally, consider those terms with $|j-l|<T$ and $|i-k|<T$. Without
loss, assume $i\leq k$. If $l\leq j$, then
%e6 ###
%
\begin{eqnarray}\label{eqn5}
&&\operatorname{Prob}(X_i=X_j, X_k=X_l) \nonumber
\\
&&\qquad = \sum_{u,v} \operatorname{Prob}(X_i=u)P^{k-i}(u,v)P^{l-k}(v,v)P^{j-l}(v,u)
\\
&&\qquad \leq \sum_u \operatorname{Prob}(X_i=u)\sum_v P^{k-i}(u,v)\frac MN P^{j-l}(v,u).\nonumber
% &\leq& \frac MN \sum_u \operatorname{Prob}(X_i=u) \sum_v P^{k-i}(u,v)P^{j-l}(v,u) .
\end{eqnarray}
\end{Case}

The sum over elements with $i\leq k< i+T$ and $l\leq j< l+T$ is upper
bounded as follows:
%e7 ###
%
\begin{eqnarray}\label{eqn:simplify-later}
 &&\sum_{i=0}^{\beta\sqrt{N}} \sum_{k=i}^{i+T} \sum
_{l=k+2T}^{\beta\sqrt{N}+2T}\sum_{j=l}^{l+T} \operatorname{Prob}(X_i=X_j, X_k=X_l) \nonumber
\\
&&\qquad \leq \frac MN \sum_{i=0}^{\beta\sqrt{N}} \sum_{l=i+2T}^{\beta
\sqrt{N}+2T} \max_u \sum_v \sum_{k\in[i,i+T)} P^{k-i}(u,v)
 \sum_{j\in[l,l+T)} P^{j-l}(v,u)\nonumber
 \\
&&\qquad \leq \frac MN \sum_{i=0}^{\beta\sqrt{N}} \sum_{l=i+2T}^{\beta
\sqrt{N}+2T} \max_u \sum_v G_T(u,v) G_T(v,u)
\\
&&\qquad \leq \frac MN \sum_{i=0}^{\beta\sqrt{N}} \sum_{l=i+2T}^{\beta
\sqrt{N}+2T} \max_u \sqrt{\sum_v G^2_T(u,v) \sum_v G^2_T(v,u)}\nonumber
\\
&&\qquad \leq \frac MN \pmatrix{\beta\sqrt{N}+2\cr2} \sqrt{A_T A_T^*} .\nonumber
\end{eqnarray}

The case when $j< l$ gives the same bound but with the observation that
$j\geq k+T$ and with $A_T$ instead of $\sqrt{A_T A_T^*}$.

Putting together these various cases we get that
\begin{eqnarray*}
E[S^2] &\leq& \pmatrix{\beta\sqrt{N}+2\cr 2}^2 \biggl(\frac
{M}{N}\biggr)^2
+ 2 \pmatrix{\beta\sqrt{N}+2\cr 2} \frac{M}{N} A_T
\\
&&{} + 2 \pmatrix{\beta\sqrt{N}+2\cr 2} \frac{M}{N} \sqrt{A_T
A_T^*} .
\end{eqnarray*}
The ${{\beta\sqrt{N}+2}\choose {2}}^2$ term is the total number of
values of $i,j,k,l$ appearing in the sum for $E[S^2]$, and hence also
an upper bound on the number of values in Cases~\ref{ca1} and~\ref{ca2}. Along with
the relation
${{\beta\sqrt{N}+2}\choose {2}}\geq\frac{\beta^2N}{2}$
this simplifies to complete the proof.
\end{pf*}

As promised earlier, we now present an example that illustrates the
need for the pre-mixing term $A_T$ in Theorem \ref{theorem:birthday}.

\begin{example} \label{ex:pre-mixing}
Consider the random walk on $\mathbb{Z}_N$ which transitions from
$u\to u+1$
with probability $1-1/\sqrt{N}$, and with probability $1/\sqrt{N}$
transitions $u\to v$ for a uniformly random choice of $v$.

Heuristically the walk proceeds as $u\to u+1$ for $\approx\sqrt{N}$
steps, then randomizes, and then proceeds as $u\to u+1$ for another
$\sqrt{N}$ steps. This effectively splits the state space into $\sqrt
{N}$ blocks of size about $\sqrt{N}$ each, so by the standard Birthday
Paradox it should require about $\sqrt{N^{1/2}}$ of these
randomizations before a collision will occur, in short, about $N^{3/4}$
steps in total.

To see the need for the pre-mixing term, observe that $T_s \approx
\sqrt{N} \log2$ while if $T=T_\infty\approx\sqrt{N} \log
(2(N-1))$, then we may take $m=1/2$ and $M=3/2$ in Theorem \ref
{theorem:birthday}. So, whether $T_s$ or $T_\infty$ are considered, it
will be insufficient to take $O(T + \sqrt{N})$ steps.
However, the number $A_T$ of collisions between two independent copies
of this walk is about $\sqrt{N}$ since once a randomization step
occurs then the two independent walks are unlikely to collide anytime soon.
Our collision time bound says that $O(N^{3/4})$ steps will suffice
which is the correct bound.

A proper analysis shows that $\frac{1-o(1)}{\sqrt{2}}N^{3/4}$ steps
are necessary to have a collision with a probability of $1/2$.
Conversely, when $T=\sqrt{N}\log^2 N$ then $m=1-o(1)$, $M=1+o(1)$ and
$A_T,A_T^*\leq\frac{1+o(1)}{2} \sqrt{N}$, so by equation (\ref
{eqn:positive_prob}), $(2+o(1))N^{3/4}$ steps are sufficient to have a
collision with a probability of at least $1/2$. Our upper bound is thus
off by at most a factor of $2\sqrt{2}\approx2.8$.
\end{example}

We finish the section with a slight sharpening of Theorem \ref
{theorem:birthday}. This will be used to improve the lead constant in our
upcoming bound on collision time for the Pollard Rho walk.

\begin{theorem}[(Improved Birthday Paradox)] \label{theorem:birthday2}
Consider a finite ergodic Markov chain with uniform stationary
distribution on a state space of size $N$. Let $T$ be such that $\frac
{m}{N} \leq P^T(u,v) \leq\frac{M}{N}$ for some $m\leq1\leq M$ and
every pair of states $u,v$.
After
\[
2c\Biggl(\sqrt{\Biggl(1+\sum_{j=1}^{2T} 3j\max_{u,v} P^j(u,v)
\Biggr)\frac NM}+T\Biggr)
\]
steps a collision occurs with a probability of at least $1-
(1-\frac{m^2}{2M^2})^c$, independent of the starting state.
\end{theorem}

\begin{pf}
We give only the steps that differ from before.

First, in equation (\ref{eqn:simplify-later}), note that the triple
sum after $\max_u$ can be re-written as
\[
\sum_{\alpha\in[0,T)}\sum_{\beta\in[0,T)} \sum_v P^{\alpha
}(u,v)P^{\beta}(v,u)
\leq\sum_{\gamma=0}^{2(T-1)} (\gamma+1)P^\gamma(u,u)  .
\]
The original quadruple sum then reduces to
\[
\frac MN\pmatrix{\beta\sqrt{N}+2\cr 2} \max_u \sum_{\gamma
=0}^{2(T-1)} (\gamma+1)P^\gamma(u,u) .
\]

For the case when $i<k$ and $j<l$ proceed similarly, then reduce, as in
Lemma~\ref{lem:quick-convergence}, to obtain the upper bound,
\begin{eqnarray*}
&&\frac MN\pmatrix{\beta\sqrt{N}+2\cr 2} \sum_{\alpha
=1}^{T-1}\sum_{\beta=1}^{T-1}\sum_v P^{\alpha}(u,v)P^{\beta}(u,v)
\\
&&\qquad \leq \frac MN\pmatrix{\beta\sqrt{N}+2\cr 2} \sum_{\gamma=1}^{T-1}
(2\gamma-1)\max_v P^\gamma(u,v) .
\end{eqnarray*}
Adding these two expressions gives an expression of at most
\[
\frac MN\pmatrix{\beta\sqrt{N}+2\cr 2}\Biggl(1+\sum_{\gamma=1}^{2T}
3\gamma\max_v P^\gamma(u,v)\Biggr) .
\]

The remaining two cases will add to the same bound, so effectively this
substitutes the expression $2(1+\max_u \sum_{\gamma=1}^{2T}
3\gamma\max_v P^\gamma(u,v))$ in place\vspace{1pt} of a $4\max\{
A_T,A_T^*\}$ in the original theorem.
\end{pf}

To simplify, note that if $\max_{u,v}P^j(u,v) \leq c+d^j$ for $c,d\in
[0,1)$, then\vspace{-2pt}
%e8 ###
%
\begin{eqnarray}\label{eqn:simplify-birthday2}
\sum_{j=1}^{2T} 3j(c+d^j)
&=& 3cT(2T+1) + 3d\frac{1-d^{2T}-2Td^{2T}(1-d)}{(1-d)^2} \nonumber
\\[-9pt]\\[-9pt]
&\leq& \bigl(1+o(1)\bigr)6cT^2 + \frac{3d}{(1-d)^2}  .\nonumber
\end{eqnarray}\vspace{-2pt}

%**************** Mixing Time ******************

%s4 ###
\section{Convergence of the Rho walk} \label{sec:mixing}

Let us now turn our attention to the Pollard Rho walk for discrete
logarithm. To apply the collision time result we will first show that
$\max_{u,v\in\mathbb{Z}_N} P^s(u,v)$ decreases quickly in $s$ so
that Lemma
\ref{lem:quick-convergence} may be used. We then find $T$ such that
$P^T(u,v)\approx1/N$ for every $u,v\in\mathbb{Z}_N$. However,
instead of
studying the Rho walk directly, most of the work will instead involve a
``block walk'' in which only a certain subset of the states visited by
the Rho walk are considered.\vspace{-2pt}

\begin{definition}
Let us refer to the three types of moves that the Pollard Rho random
walk makes,
namely $(u, u+1), (u,u+x)$ and $(u, 2u)$ as moves of Type~1, Type 2 and
Type 3,
respectively. In general, let the random walk be denoted by $Y_0, Y_1,
Y_2, \ldots$  with
$Y_t$ indicating the position of the walk (modulo $N$) at time $t\ge
0$. Let $T_1$ be the first
time that the walk makes a move of Type 3.
Let $b_1 = Y_{T_1 - 1} - Y_{T_0}$ (i.e., the ground covered, modulo
$N$, only using
consecutive moves of Types 1 and 2). More generally, let $T_i$ be the
first time, since $T_{i-1}$, that a move of Type 3 happens, and set
$b_i = Y_{T_i - 1} - Y_{T_{i-1}}$. Then the \textit{block walk} $\mathsf
{B}$ is
the walk $X_s = Y_{T_s} = 2^s Y_{T_0} + 2\sum_{i=1}^s 2^{s-i}b_i$.
\end{definition}\vspace{-2pt}

By combining our Birthday Paradox for Markov chains with several lemmas
to be shown in this section, we obtain the main result of the paper:\vspace{-2pt}

\begin{theorem} \label{theorem:rho-collision}
For every choice of starting state, the expected number of steps
required for the Pollard Rho algorithm for discrete logarithm on a
group $G$ to have a collision is at most\vspace{-2pt}
\[
\bigl(1+o(1)\bigr) 12\sqrt{19} \sqrt{|G|} < \bigl(1+o(1)\bigr) 52.5 \sqrt{|G|} .
\]\vspace{-2pt}
\end{theorem}\vspace{-2pt}

In order to prove this it is necessary to show that $\mathsf{B}^s(u,v)$
decreases quickly for the block walk.\vspace{-2pt}

\begin{lem} \label{lem:rho-quick}
If $s\leq\lfloor\log_2 N\rfloor$, then for every $u,v\in\mathbb{Z}_N$
the block walk satisfies\vspace{-2pt}
\[
\mathsf{B}^s(u,v) \leq(2/3)^s .\vspace{-2pt}
\]
If $s>\lfloor\log_2 N\rfloor$ then $\displaystyle\mathsf
{B}^s(u,v)\leq
\frac{3/2}{N^{\log_2 3-1}} \leq\frac{3/2}{\sqrt{N}}$.
\end{lem}\vspace{-3pt}

A bound on the asymptotic rate of convergence is also required:\vspace{-3pt}

\begin{theorem} \label{theorem:asymptotics}
If $s\geq\lceil m\log\frac{2(m-1)}{\varepsilon}\rceil$
where $m=\lceil\log_2 N\rceil$, then for every $u,v\in
\mathbb{Z}_N$ the block walk satisfies\vspace{-2pt}
\[
\frac{1-\varepsilon}{N}\leq\mathsf{B}^{2s}(u,v)\leq\frac{1+\varepsilon
}{N} .
\]
\end{theorem}\vspace{-3pt}

This is all that is needed to prove the main result:\vspace{-2pt}

\begin{pf*}{Proof of Theorem \ref{theorem:rho-collision}}
The proof will use Theorem \ref{theorem:birthday2} because this gives a
somewhat sharper bound. Alternatively, Theorem \ref{theorem:birthday} and
Lemma \ref{lem:quick-convergence} can be applied nearly identically to
get the slightly weaker $(1+o(1))72\sqrt{|G|}$.

First consider steps of the block walk. Lemma \ref{lem:rho-quick}
implies that
$\mathsf{B}^s(u,v) \le\frac{3/2}{\sqrt{N}} + {(\frac{2}{3})}^s$,
for $s\ge
0$ and for all $u,v$.
Hence, by equation (\ref{eqn:simplify-birthday2}), if $T=o(\sqrt
[4]{N})$, then $1+\sum_{j=1}^{2T} 3j \mathsf{B}^j(u,v) \leq
19+o(1)$. By
Theorem \ref{theorem:asymptotics}, $M\leq1+\varepsilon$ and $m\geq
1-\varepsilon$ after $2(\log_2 N)(\log\log N + \log\frac{3}{\varepsilon
})$ steps. Hence, if $\varepsilon=1/N^2$, then $T=(4+o(1)) (\log_2
N)^2=o(\sqrt[4]{N})$, and $m=1-o(1/N)$, and $M=1+o(1/N)$. Plugging
this into Theorem \ref{theorem:birthday2}, a collision fails to occur in\vspace{-2pt}
\[
k\Biggl(2\sqrt{\Biggl( 1+\sum_{j=1}^{2T} 3j \max_{u,v} \mathsf{B}
^j(u,v)\Biggr)\frac NM}+2T\Biggr)
= \bigl(1+o(1)\bigr) 2\sqrt{19} k\sqrt{N}\vspace{-2pt}
\]
%
% Lemma \ref{lem:quick-convergence},
%%we may set $c=\frac{3/2}{\sqrt{N}}$ and $d=2/3$. Hence,
% if $T=o(\sqrt[4]{N})$ then $A_T,A_T^*\leq\frac{2+o(1)}{(1 -
%2/3)^2}=(1+o(1)) 18$. By Remark \ref{rmk:L2-to-pointwise}, after $2(
%1+\varepsilon$ and $m\geq1-\varepsilon$. Hence, if we set $\varepsilon=1/N^2$
%then $T=(4+o(1)) (\log_2 N)^2$ with $m=1-o(1/N)$ and $M=1+o(1/N)$.
%Plugging this into equation \eqref{eqn:positive_prob} and the
%subsequent few lines, a collision fails to occur in
%$$
%k(2\sqrt{\frac{2N}{M}\max\{A_T, A_T^*\}}+2T)
% = (1+o(1)) 12k\sqrt{N}
%$$
steps with a probability of, at most, $(1-\delta)^k$ where $\delta
=m^2/2M^2=(1-o(1))/2$.\vspace{-2pt}

Now return to the Rho walk. Recall that $T_i$ denotes the number of Rho
steps required for $i$ block steps. The difference $T_{i+1}-T_i$ is an
i.i.d. random variable with the same distribution as $T_1-T_0$. Hence,
if $i\geq j$, then $E[T_i-T_j] = (i-j) E[T_1-T_0]=3(i-j)$. In
particular, if we let $r=(1+o(1)) 2\sqrt{19 N}$, let $R$ denote the
number of Rho steps before a collision, and let $B$ denote the number
of block steps before a collision, then\vspace{-2pt}
\begin{eqnarray*}
E[R] &\leq& \sum_{k=0}^{\infty} \operatorname{Pr}[B>kr] E\bigl[T_{(k+1)r}-T_{kr} \mid B>kr\bigr]
\\[-3pt]
&=& \sum_{k=0}^{\infty} \operatorname{Pr}[B>kr] E\bigl[T_{(k+1)r}-T_{kr}\bigr]
\\[-3pt]
&\leq& \sum_{k=0}^{\infty} \biggl(\frac{1+o(1)}{2}\biggr)^k 3r
= \bigl(1+o(1)\bigr) 12\sqrt{19} \sqrt{N} .\vspace*{-8pt}
\end{eqnarray*}
\upqed\end{pf*}

\begin{pf*}{Proof of Lemma \ref{lem:rho-quick}} We start with a
weaker, but somewhat more intuitive, proof of a bound on $\mathsf{B}^s(u,v)$
and then improve it to obtain the result of the lemma. The key idea
here will be to separate out a portion of the Markov chain which is
tree-like with some large depth $L$, namely the moves induced solely by
$b_i=0$ and $b_i=1$ moves. Because of the high depth of the tree, the
walk spreads out for the first $L$ steps, and hence the probability of
being at a vertex also decreases quickly.

Let $S=\{i\in[1,\ldots, s] \dvtx b_i\in\{0,1\}\}$ and $z=\sum_{i\notin
S} 2^{s-i}b_i$ be random variables whose values are determined by the
first $T_s$ steps of the random walk. Then $Y_{T_s}=2^s Y_{T_0} + 2z +
2 \sum_{i\in S} 2^{s-i}b_i$. Hence, choosing $Y_{T_0}=u$,
$Y_{T_s} = v$, we may write
\begin{eqnarray*}
&&\mathsf{B}^s(u,v)
\\
&&\qquad = \sum_{S} \operatorname{Prob}(S)
\sum_{z\in\mathbb{Z}_N} \operatorname{Prob}(z\mid S) \operatorname{Prob}\Biggl(\sum_{i\in S} 2^{s-i}b_i
= \frac v2-2^{s-1} u-z  \Bigm| z, S\Biggr)
\\
&&\qquad \leq \sum_{S} \operatorname{Prob}(S)  \max_{w\in\mathbb{Z}_N} \operatorname{Prob}\biggl(\sum
_{i\in
S} 2^{s-i}b_i=w  \Bigm|  S\biggr) ,
\end{eqnarray*}
and so for a fixed choice of $S$, we can ignore what happens on $S^c$.

Each $w\in[0,\ldots, N-1]$ has a unique binary expansion, and so if
$s\leq\lfloor\log_2 N\rfloor$, then modulo $N$ each $w$ can still
be written in, at most, one way as an $s$ bit string. For the block
walk, $\operatorname{Prob}(b_i=0)\geq1/3$ and $\operatorname{Prob}(b_i=1) \geq1/9$, and so $\max\{
\operatorname{Prob}(b_i=0\mid i\in S), \operatorname{Prob}(b_i=1\mid i\in S)\}\leq\frac89$. It
follows that
%e9 ###
%
\begin{equation} \label{eqn:maxprob}
\max_{w\in\mathbb{Z}_N} \operatorname{Prob}\biggl(\sum_{i\in S} 2^{s-i}b_i=w
\Bigm|
S\biggr) \leq(8/9)^{|S|}
\end{equation}
using independence of the $b_i$s.
Hence,
\begin{eqnarray*}
\mathsf{B}^s(u,v) &\leq& \sum_{S} \operatorname{Prob}(S) (8/9)^{|S|}
= \sum_{r=0}^s \operatorname{Prob}(|S|=r) (8/9)^r
\\
& \leq& \sum_{r=0}^s \pmatrix{s\cr r}\biggl(\frac49\biggr)^r
\biggl(1-\frac49\biggr)^{s-r} \biggl(\frac89\biggr)^r
= \biggl(\frac49\frac89+\frac59\biggr)^s = \biggl(\frac
{77}{81}\biggr)^s .
\end{eqnarray*}
The second inequality was because $(8/9)^{|S|}$ is decreasing in $|S|$
and so underestimating $|S|$ by assuming $\operatorname{Prob}(i\in S)=4/9$ will only
increase the upper bound on $\mathsf{B}^s(u,v)$.

In order to improve on this, we will shortly redefine $S$ (namely,
events $\{i \in S\}, \{i \notin S\}$) and auxiliary variables $c_i$,
using the steps of the Rho walk. Also note that the block walk is
induced by a Rho walk, so we may assume that the $b_i$ were constructed
by a series of steps of the Rho walk. With probability $1/4$ set $i\in
S$ and $c_i=0$, otherwise if the first step is of Type 1, then set
$i\in S$ and $c_i=1$ while if the first step is of Type 3 then put
$i\notin S$ and $c_i=0$, and finally if the first step is of Type~2,
then again repeat the above decision making process, using the
subsequent steps of the walk. Note that the above construction can be
summarized as consisting
of one of four equally likely outcomes (at each time) where the last
three outcomes
depend on the type of the step that the Rho walk takes; indeed, each of these
three outcomes happens with probability $\frac{3}{4} \times\frac
{1}{3} = 1/4$; finally,
a Type 2 step forces us to reiterate the four-way decision-making process.

In summary, $\operatorname{Pr}(i\in S)=\sum_{l=0}^{\infty} (1/4)^l (1/2)=2/3$.
Also observe that $\operatorname{Pr}(c_i=0\mid i\in S) = \operatorname{Pr}(c_i=1\mid i\in S)$ and that
$\operatorname{Pr}(b_i-c_i=x\mid i\in S, c_i=0)=\operatorname{Pr}(b_i-c_i=x\mid i\in S, c_i=1)$.
Hence the steps done earlier (leading to the weaker bound) carry
through with $z=\sum_{i} 2^{s-i}(b_i-c_i)$ and with $\sum_{i\in
S}2^{s-i} b_i$ replaced by $\sum_{i\in S}2^{s-i} c_i$.
In (\ref{eqn:maxprob}) replace $(8/9)^{|S|}$ by $(1/2)^{|S|}$, and in
showing the final upper bound on $\mathsf{B}^s(u,v)$, replace $4/9$ by $2/3$.
This leads to the bound $\mathsf{B}^s(u,v)\leq(2/3)^s$.

Finally, when $s>\lfloor\log_2 N\rfloor$, simply apply the preceding
argument to $S'=S\cap[1,\ldots,\lfloor\log_2 N\rfloor]$.
Alternately, note that when $s\geq\lfloor\log_2 N\rfloor$, then
\[
\mathsf{B}^s(u,v) \leq\max_w \mathsf{B}^{\lfloor\log_2 N\rfloor
}(u,w) \leq
(2/3)^{\log_2 N - 1}
\]
for every doubly-stochastic Markov chain $\mathsf{B}$.
\end{pf*}

%As was demonstrated in the above proofs, in order to use the Birthday
%Paradox on the Rho walk it suffices to show a mixing time bound of
%$T=O(\sqrt[4]{N})$ (to guarantee that $A_T,A_T^*=O(1)$).
In \cite{MV061,KMT071} sufficiently strong bounds on the asymptotics
of $\mathsf{B}^{2s}(u,v)$ are shown in several ways, including the use of
characters and quadratic forms, canonical paths and Fourier analysis.
We give here the Fourier approach, as it establishes the sharpest
mixing bounds.
To bound mixing time of the block walk, it suffices to show that for
large enough $s$, the distribution $\nu_s $ of
\[
Z_s = 2^{s-1} b_1 + 2^{s-2} b_2 + \cdots+ b_{s}
\]
is close to the uniform distribution $U=1/N$ because then the
distribution of $X_s=2^s Y_{T_0}+2Z_s$ will be close to uniform as
well. More precisely, convergence in chi-square distance will be shown that

\begin{lem} \label{lem:fourier}
If $\nu_s (j) = \mathsf{\operatorname{Pr}}[ Z_s = j]$, $\xi= 1- \frac{4-\sqrt{10}}{9}$,
and $m$
satisfies $2^{m-1} < N < 2^m$, then
\[
N \sum_{j=0}^{N-1} \bigl(\nu_s (j) - U(j)\bigr)^2
\leq2 \bigl( \bigl(1+\xi^{2\lfloor s/m \rfloor}\bigr)^{m-1} -1 \bigr).
\]
\end{lem}

\begin{pf*}{Proof of Theorem \ref{theorem:asymptotics}}
By Cauchy--Schwarz,
%e10 ###
%
\begin{eqnarray}\label{eqn:cauchy-schwartz}
&&\biggl| \frac{\mathsf{B}^{2s}(u,v)-U(v)}{U(v)} \biggr|^2\nonumber
\\
&&\qquad= \biggl| \frac{\sum_w (\mathsf{B}^s(u,w)-U(w))
(\mathsf{B}
^s(w,v)-U(v))}{U(v)} \biggr|^2 \nonumber
\\[-8pt]\\[-8pt]
&&\qquad = \biggl| \sum_w U(w) \biggl(\frac{\mathsf
{B}^s(u,w)}{U(w)}-1
\biggr)\biggl(\frac{\mathsf{B}^{*s}(v,w)}{U(w)}-1\biggr)\biggr|^2\nonumber
\\
&&\qquad \leq \sum_w U(w) \biggl|\frac{\mathsf{B}^s(u,w)}{U(w)}-1
\biggr|^2 \sum
_x U(x)\biggl|\frac{\mathsf{B}^{*s}(v,x)}{U(x)}-1\biggr|^2.\nonumber
\end{eqnarray}
Lemma \ref{lem:fourier} bounds the first sum of
(\ref{eqn:cauchy-schwartz}). The second sum is the same quantity but for the
time-reversed walk $\mathsf{B}^*(y,x)=\mathsf{B}(x,y)$. To examine
the reversed walk
let $b_i^*$ denote the sum of steps taken by $\mathsf{B}^*$ between the
$(i-1)$st and $i$th time that a $u\to u/2$ transition is chosen (i.e.,
consider block steps for the reversed walk), and let $Z_s^*=2^{-s+1}
b_1^*+\cdots+b_s^*$. If we define $b_i=-b_i^*$, then the $b_i$ are
independent random variables from the same distribution as the blocks
of $\mathsf{B}$, and so
\begin{eqnarray*}
\mathsf{Pr}[-2^{s-1}Z_s^*=j] &=& \mathsf{Pr}[b_1+2b_2+\cdots
+2^{s-1}b_s=j] \\
&=& \mathsf{Pr}[Z_s=j] .
\end{eqnarray*}
Lemma \ref{lem:fourier} thus bounds the second sum of
(\ref{eqn:cauchy-schwartz}) as well, and the theorem\break
follows.
%
%It follows from \eqref{eqn:L2} that
%$$
%| 1-\frac{\\operatorname{Pr}[Z_{2s}=j]}{U(j)} |
% \leq2 ( (1+\xi^{2\lfloor s/m \rfloor})^{m-1} -1 ) ,
%$$
%and so after $2s\approx m \log_{\xi}\frac{\varepsilon/2}{m-1}\leq2m\log
\end{pf*}

Before proving Lemma \ref{lem:fourier} let us review the standard
Fourier transform and the Plancherel identity. For any complex-valued
function $f$ on $\mathbb{Z}_N$ and $\omega= e^{2\pi i/N}$,
recall that the Fourier transform $\hat{f}\dvtx \mathbb{Z}_N\rightarrow C$
is given by
$\hat{f} (\ell) = \sum_{j=0}^{N-1} \omega^{\ell j }
f(j)$, and the
Plancherel identity asserts that
\[
N \sum_{j=0}^{N-1} |f(j)|^2= \sum_{j=0}^{N-1} |\hat{f}(j)|^2 .
\]

For the distribution $\mu$ of a $\mathbb{Z}_N$-valued random variable
$X$, its
Fourier transform is
\[
\hat{\mu} (\ell)
= \sum_{j=0}^{N-1} \omega^{\ell j } \mu(j)
= E[ \omega^{\ell X}].
\]
Thus, given distributions $\mu_{_1}, \mu_{_2}$ of two independent
random variables $Y_1, Y_2$, the distribution $\nu$ of
$X:=Y_1+Y_2$ has the Fourier transform $\hat{\nu}=\hat{\mu
}_{_1}\hat{\mu}_{_2}$, since

%e11 ###
%
\begin{eqnarray}\label{eqn:convolution}
\hat{\nu} (\ell)
&=& E[ \omega^{\ell X}]
= E\bigl[ \omega^{\ell(Y_1+Y_2)}\bigr] \nonumber
\\[-8pt]\\[-8pt]
&=& E[ \omega^{\ell Y_1}]E[ \omega^{\ell Y_2}]
= \hat{\mu}_{_1} (\ell)\hat{\mu}_{_2}(\ell).\nonumber
\end{eqnarray}
Generally, the distribution $\nu$ of $X:=Y_1 +\cdots+ Y_s$
with independent $Y_i$s has the Fourier transform
$\hat{\nu}=\prod_{r=1}^s \hat{\mu}_{_r}$.
Moreover, for the uniform distribution $U$, it is easy to check that
\[
\hat{U} (\ell)
= \cases{ {1},&\quad \mbox{if } $\ell=0$, \cr {0}, &\quad \mbox{otherwise.}
}
\]
As the random variables $2^{r}b_{s-r}$s are independent,
$\hat{\nu}_s = \prod_{r=0}^{s-1} \hat{\mu}_r $ where $\mu_r$ are
the distributions of $2^{r}b_{s-r}$. The linearity of the Fourier
transform and
$\hat{\nu}_s (0) = E[1]=1$ yield
\[
\widehat{\nu_s - U} (\ell)
= \hat{\nu}_s (\ell)- \hat{U}(\ell)
= \cases{ {0}, &\quad \mbox{if } $\ell=0$, \cr {\displaystyle\prod_{r=0}^{s-1}
\hat{\mu}_r(\ell)},&\quad \mbox{otherwise.}}
\]

\begin{pf*}{Proof of Lemma \ref{lem:fourier}}
By Plancherel's identity, it is enough to show that
\[
\sum_{\ell=1}^{N-1} \Biggl|\prod_{r=0}^{s-1} \hat{\mu}_r
(\ell) \Biggr|^2 \leq2 \bigl( \bigl(1+\xi^{2\lfloor s/m \rfloor}\bigr)^{m-1}
-1 \bigr).
\]

Let $A_r$ be the event that $b_{s-r} = 0 ~{\rm or} ~1$. Then,
\begin{eqnarray*}
\hat{\mu}_r(\ell)
&=& E [ \omega^{\ell2^{r} b_{s-r}}]
\\
&=& \mathsf{Pr}[ b_{s-r}=0] + \mathsf{Pr}[b_{s-r}=1] \omega^{\ell 2^{r} }
\\
&& {} + \mathsf{Pr}[ \bar{A}_r ]E [ \omega^{\ell2^{r} b_{s-r}}| \bar{A}_r],
\end{eqnarray*}
and, for $x:=\mathsf{\operatorname{Pr}}[b_{s-r}=0]$ and $y:= \mathsf{\operatorname{Pr}}[b_{s-r}=1]$,
\begin{eqnarray*}
|\hat{\mu}_r(\ell)|
&\leq& | x + y \omega^{\ell2^{r} }| + (1-x-y)| E [ \omega^{\ell2^{r}
b_{s-r}}| \bar{A}_r]|
\\
&\leq& | x + y \omega^{\ell2^{r} }| + 1-x-y .
\end{eqnarray*}
Notice that
\begin{eqnarray*}
| x + y \omega^{\ell2^{r} }|^2
&=& \biggl( x+ y \cos\frac{2\pi\ell2^r}{N} \biggr)^2+ y^2 \sin^2
\frac{2\pi\ell2^r}{N}
\\
&=& x^2 + y^2 + 2 xy \cos\frac{2\pi\ell2^r}{N}.
\end{eqnarray*}
If $\cos\frac{2\pi\ell2^r}{N}\leq0$, then
\begin{eqnarray*}
|\hat{\mu}_r (\ell)|
&\leq& (x^2 +y^2)^{1/2}+1-x-y
\\
&=& 1-\bigl(x+y-(x^2 +y^2)^{1/2}\bigr).
\end{eqnarray*}
Since $x= \mathsf{\operatorname{Pr}}[b_ {s-r}=0] \geq1/3$, and $y= \mathsf
{\operatorname{Pr}}[b_{s-r}=1]\geq1/9$,
it is easy to see that $x+y-(x^2 +y^2)^{1/2}$ has its minimum when
$x=1/3$ and $y=1/9$.
(For both partial derivatives are positive.) Hence,
\[
|\hat{ \mu}_r (\ell)| \leq\xi=1- \frac{4-\sqrt{10}}{9}
\qquad\mbox{provided } \cos\frac{2\pi\ell2^r}{N} \leq0.
\]
If $ \cos\mbox{$\frac{2\pi\ell2^r}{N}$}> 0$, we use the trivial bound
$\hat{\mu}_r (\ell) = E [ \omega^{\ell2^{r} b_{s-r}}] \leq1$.

For $\ell=1,\ldots, N-1$, let $\phi_s (\ell) $ be the number of
$r=0,\ldots, s-1$ such that
$\cos\frac{2\pi\ell2^r}{N}\leq0$. Then
\begin{eqnarray}\label{phibd}
\prod_{r=0}^{s-1} | \hat{\mu}_r
(\ell)| \leq\xi^{\phi_s (\ell)}.
\end{eqnarray}
To estimate $\phi_s (\ell)$, we consider the binary expansion of
\[
\ell/ N = 0. \alpha_{_{\ell,1}} \alpha_{_{\ell,2}} \cdots\alpha
_{_{\ell
,s}} \cdots,
\]
$\alpha_{_{\ell,r}} \in\{0,1\} $ with $\alpha_{_{\ell,r}} =0$
infinitely often. Hence, $\ell/N = \sum_{r=1}^\infty2^{-r}
\alpha_{_{\ell,r}}$. The fractional part of $\ell2^{r} / N$ may be
written as
\[
\{ \ell2^{r} / N\} =0 .\alpha_{_{\ell, r+1}}\alpha_{_{\ell, r+2}}
\cdots
\alpha_{_{\ell,s}} \cdots.
\]
Notice that $\cos\frac{2\pi\ell2^{r}}{ N} \leq0$ if the
fractional part of $\ell2^{r} / N$ is (inclusively) between
$1/4$ and $3/4$ which follows if $\alpha_{_{r+1}} \not=
\alpha_{_{r+2}}$. Thus $\phi_s (\ell) $ is at least as large as the
number of
alterations in the sequence $(\alpha_{_{\ell,1}}, \alpha_{_{\ell,2}},
\ldots, \alpha_{_{\ell, s+1}})$.

We now take $m$ such that $2^{m-1} < N < 2^m$. Observe that, for
$\ell=1, \ldots , N-1$, the subsequences $\alpha(\ell) :=
(\alpha_{_{\ell,1}}, \alpha_{_{\ell,2}}, \ldots , \alpha_{_{\ell,
m}})$ of
length $m$ are pairwise distinct: If $\alpha(\ell)=\alpha(\ell')$ for
some $\ell< \ell'$ then $\frac{ \ell'-\ell}{N}$ is less than
$\sum_{r\geq m+1} 2^{-r} \leq2^{-m}$ which is impossible as $N <
2^m$. Similarly, for fixed $r$ and $\ell=1,\ldots , N-1$, all
subsequences $\alpha(\ell;r ) :=
(\alpha_{_{\ell,r+1}}, \alpha_{_{\ell,r+2}}, \ldots , \alpha_{_{\ell, r+m}})$
are pairwise distinct. In
particular, for fixed $r$ with $r=0, \ldots, \lfloor s/m \rfloor-1$,
all subsequences $\alpha(\ell;rm) $, $\ell=1,\ldots , N-1$, are pairwise
distinct. Since the fractional part $\{ \frac{2^{rm} \ell}{N} \}=
0.\alpha_{_{\ell,rm+1}} \alpha_{_{\ell,rm+2}} \cdots$ must be the
same as
$\frac{ \ell' }{N}$ for some $\ell'$ in the range $1\leq\ell'
\leq N-1$, there is a unique permutation $\sigma_r$ of $1, \ldots, N-1$
such that $\alpha(\ell;rm)= \alpha( \sigma_r (\ell))$. Writing
$|\alpha(
\sigma_r (\ell))|_{_A} $ for the number of alternations in $\alpha(
\sigma_r (\ell))$, we have
\[
\phi_s (\ell) \geq\sum_{r=0}^{\lfloor s/m \rfloor-1} |\alpha(
\sigma_r
(\ell))|_{_A},
\]
where $\sigma_{_0}$ is the identity. Therefore, (\ref{phibd}) gives
\[
\sum_{\ell=1}^{N-1} \Biggl|\prod_{r=0}^{s-1} \hat{\mu}_r
(\ell) \Biggr|^2 \leq\sum_{\ell=1}^{N-1}
%%%\prod_{r=0}^{\lfloor s \rfloor-1}
  \xi^{2 \sum_{r=0}^{\lfloor s/m \rfloor-1} |\alpha( \sigma_r
(\ell))|_{_A}}.
\]

Using
\begin{eqnarray*}
&&\xi^{x+y}+ \xi^{x'+y'}
 \\
&&\qquad \leq \xi^{\min\{x,x'\} +\min\{y,y'\}} +
\xi^{\max\{x,x'\} +\max\{y,y'\}}
\end{eqnarray*}
inductively, the above upper bound may be maximized when all $\sigma_r$s
are the identity, that is,
\[
\sum_{\ell=1}^{N-1} \Biggl|\prod_{r=0}^{s-1} \hat{\mu}_r
(\ell) \Biggr|^2 \leq\sum_{\ell=1}^{N-1}
%%% \prod_{r=0}^{\lfloor s \rfloor-1}
  \xi^{2\lfloor s/m \rfloor|\alpha( \ell)|_{_A}}.
\]
Note that $1/N \leq\ell/N \leq1- 1/N$ implies that $\alpha(\ell)$
is neither $(0,\ldots, 0)$ nor $(1, \ldots, 1)$ (both are of length
$m$). This means that all $\alpha(\ell)$ have at least one
alternation. Since $\alpha(\ell)$s are pairwise distinct,
\[
\sum_{\ell=1}^{N-1}   \xi^{2\lfloor s/m \rfloor|\alpha(
\ell)|_{_A}} \leq\sum_{\alpha: |\alpha|_{_A} >0} \xi^{2 \lfloor s/m
\rfloor|\alpha|_{_A}},
\]
where the sum is taken over all sequences $\alpha\in\{ 0, 1\}^m$
with $ |\alpha|_{_A} >0$.

Let $H(z) $ be the number of $\alpha$s with exactly $z$ alterations.
Then
\[
H(z) =2 \pmatrix{m-1\cr z} ,
\]
and hence
\begin{eqnarray*}
\sum_{\alpha: |\alpha|_{_A} >0} \xi^{2 \lfloor s/m \rfloor|\alpha|_{_A}}
&=& 2 \sum_{z=1}^{m-1} \pmatrix{m-1\cr z} \xi^{2 \lfloor s/m \rfloor z}
\\
&=& 2 \bigl(\bigl(1+\xi^{2\lfloor s/m \rfloor}\bigr)^{m-1} -1 \bigr).
\end{eqnarray*}
\upqed\end{pf*}

\begin{rmk}
For the reader interested in applying these methods to show a
Birthday-type result for other problems, it is worth noting that a
Fourier approach can also be used to show that $\mathsf{B}^s(u,v)$ decreases
quickly, and so $A_T,A_T^*=O(1)$.

For the distribution $\nu_s$ of $X_s$ the Plancherel identity gives
\begin{eqnarray*}
\max_v \mathsf{\operatorname{Pr}}[ X_s = v]
&=& \max_v \nu_s (v)^2 \leq\sum_{w=0}^{N-1} \nu_s (w)^2
\\
&=& \frac{1}{N}\sum_{\ell=0}^{N-1} |\hat{\nu}_{s} (\ell)|^2
= \frac{1}{N}\sum_{\ell=0}^{N-1} \Biggl|\prod_{r=0}^{s-1} \hat{\mu}_r
(\ell) \Biggr|^2 .
\end{eqnarray*}

For $\ell=0, 1,\ldots, N-1$, let $\phi_s(\ell) $ be the number of
$r=0,\ldots, s-1$ such that $\cos\frac{2\pi\ell2^r}{N}\leq
0$. Then
\[
\prod_{r=0}^{s-1} | \hat{\mu}_r (\ell)|
\leq\xi^{\phi_s (\ell)}.
\]
Take $m$ such that $2^{m-1} < N < 2^m$. Then, for $s \leq m-1$
and any (fixed) binary sequence $\alpha_1, \ldots, \alpha_{s} $ (that
is, $\alpha_j \in\{0,1\}$), there are at most $\lceil2^{-s} N
\rceil$ $\ell$s such that the binary expansion of $\ell/N$ up to
$s$ digits is $0.\alpha_1 \cdots \alpha_{s} $. Since there are at
most $2e^{-\Omega(s)} 2^s$ binary sequences with fewer than
$(s-1)/3$ alterations,
\[
\prod_{r=0}^{s-1} | \hat{\mu}_r (\ell)| = 2e^{-\Omega(s)}
\]
except for at most $2e^{-\Omega(s)} 2^s \lceil2^{-s} N
\rceil=2e^{-\Omega(s)} N $ values of $\ell$. Using a trivial bound
$\prod_{r=0}^{s-1} | \hat{\mu}_r (\ell)|\leq1$ for such $\ell$'s,
we have
\[
\max_v \mathsf{\operatorname{Pr}}[ X_s = v]
= 2e^{-\Omega(s)} + 2e^{-\Omega(s)}
= 2e^{-\Omega(s)}.
\]
If $s>m-1$, then $\prod_{r=0}^{s-1} | \hat{\mu}_r (\ell)|\leq
\prod_{r=0}^{m-2} |\hat{\mu}_r(\ell)| $ implies that
\[
\max_v \mathsf{\operatorname{Pr}}[ X_s = v] =2e^{-\Omega(m-1)}= O\bigl(N^{-\Omega(1)}\bigr).
\]
\end{rmk}

One might expect that the correct order of the mixing time of the Block
walk $X_s$ is indeed $\Theta(\log p \log\log p)$. This is in fact the
case, at least for certain values of $p$ and $x$.

\begin{theorem} \label{theorem:lower-bound}
If $p=2^t - 1$ and $x=p-1$, then the block walk has mixing time $\tau
(1/4)=\Theta(\log p \log\log p)$.
\end{theorem}

\begin{pf}
The upper bound on mixing time, $O(\log p\log\log p)$, was shown in
Theorem \ref{theorem:asymptotics} via a Fourier argument.

The proof of a lower bound of $\Omega(\log p\log\log p)$ is modeled
after an arguments of Hildebrand \cite{Hil061} which in turn closely
follows a proof of Chung, Diaconis and Graham \cite{CDG871}. The
basic idea is by now fairly standard: choose a function and show that
its expectation under the stationary distribution and under the
$n$-step distribution $P^n$ are far apart with sufficiently small
variance to conclude that the two distributions ($P^n$ and $\pi$) must
differ significantly. Theorem \ref{theorem:lower-bound} is not used in
main results of this paper, and the proof is fairly long, and so it is
left for the \hyperref[app]{Appendix}.
\end{pf}

%********************** Parallel Rho ******************************

%s5 ###
\section{Distinguished point methods} \label{sec:parallel}

The Rho algorithm can be parallelized to~$J$ processors via the
Distinguished Points method of van Oorschot and Wiener \cite{VW99}. To
do this, start with a global choice of (random) partition $S_1\amalg
S_2\amalg S_3$ (i.e., a~common iterating function $F$), and choose $J$
initial values $\{y_0^j\}_{j=1}^J$ from $\mathbb{Z}_N$, one per processor.
Then run the Rho walk on processor $j$ starting from initial state
$g^{(y_0^j)}$, until a collision occurs between either two
walks or a walk and itself. To detect a collision let $\varphi\dvtx G\to
\{0,1\}$ be an easily computed hash function with support $\{x\in G\dvtx
\varphi(x)=1\}$ to be called the \textit{distinguished points}. Each time
a distinguished point is reached by a processor, it is sent to a
central repository and compared against previously received states.
Once a distinguished point is reached twice then a collision has
occurred, and the discrete logarithm can likely be found while,
conversely, once a collision occurs, the collision will be detected the
next time a distinguished point is reached.

The proofs in previous sections immediately imply a factor of $J$
speed up when parallelizing. To see this, suppose the initial values
$\{y_0^j\}_{j=1}^J$ are chosen uniformly at random. Run a Rho walk for
some $\mathcal{T}$ steps per processor, and then define $\{X_i\}$ by
starting with the Rho walk of processor $\#1$, then appending that from
processor $\#2$, etc.; that is, if $Y_i^j$ denotes the $i$th state
of copy $j$ of the walk, for $i\in\{0,1,\ldots,{\mathcal T}\}$ and $j\in
\{0,1,2,\ldots,J-1\}$, then\vspace{1pt} $X_i=Y_{i \operatorname{mod} ({\mathcal T}+1)}^{i \operatorname{div}
({\mathcal T}+1)}$ for $i\in\{0,1,\ldots,J(\mathcal{T}+1)-1\}$. This is a
time-dependent random walk which follows the Rho walk, except at
multiples of time $\mathcal{T}+1$ where it instead jumps to a uniformly
random state. Since our proofs involve pessimistic estimates on the
distance of a distribution from uniform, and these jumps result in
uniform samples, then they can only improve the result. Hence this
effectively leads to a Rho walk with $J(\mathcal{T}+1)-1$ steps, and a
factor~$J$ speed-up per processor is achieved. If the initial values
were not uniform then discard the first $O(\log^2 N)$ steps per
processor and treat the next state as the initial value which by
Theorem \ref{theorem:asymptotics} will give a nearly uniform start state.

\begin{appendix}
\section*{Appendix}\label{app}

\begin{pf*}{Proof of Theorem \ref{theorem:lower-bound}}
Our approach follows that taken in Section 4 of Hildebrand \cite{Hil061}, ``A~proof of Case
2.''
Recall that we will show that $E_U(f)$ and $E_{P^n}(f)$ are far apart
for some function $f$ with sufficiently small variance to conclude that
the two distributions ($P^n$ and the uniform distribution $U$ on
$\mathbb{Z}
_p$) must differ significantly.

More precisely let $P_n$ be the distribution of the block walk on
$\mathbb{Z}
_p$ starting at state $u=0$ and proceeding for $n$ steps.
For some $\alpha>0$, let
\[
A = \bigl\{y\dvtx |f(y)-E_U(f)| \geq\alpha\sqrt
{\operatorname{Var}
_U(f)}\bigr\}.
\]
By Chebyshev's inequality, $U(A) \leq1/\alpha^2$. Also, for some
$\beta>0$ let
\[
B = \bigl\{y\dvtx |f(y)-E_{P_n}(f)| \geq\beta\sqrt
{\operatorname{Var}
_{P_n}(f)}\bigr\} .
\]
By Chebyshev's inequality, $P_n(B) \leq1/\beta^2$.
If $A^c\cap B^c=\varnothing$, then $P_n(A^c)\leq P_n(B)\leq1/\beta^2$,
and so
\[
\min_{v\in A^c} \frac{P_n(v)}{U(v)} \leq\frac{P_n(A^c)}{U(A^c)}
\leq\frac{1/\beta^2}{1-1/\alpha^2} .
\]
If $\sqrt{\operatorname{Var}_U(f)},\sqrt{\operatorname
{Var}_{P_n}(f)}= o(|E_U(f)-E_{P_n}(f)|)$
for a sequence $n(p)=\break\Omega(\log p \log\log p)$, then as $p\to
\infty$ it is possible to choose $\alpha,\beta\to\infty$ such that
$A^c\cap B^c=0$. The theorem then follows as $\min_{v\in\Omega}
\frac{P_n(v)}{U(v)} \stackrel{p\to\infty}{\longrightarrow}0$.

The ``separating'' function $f\dvtx \mathbb{Z}_p\to\mathbb{C}$ to be used
here is
\[
f(k) = \sum_{j=0}^{t-1} q^{k2^j}\qquad\mbox{where } q=e^{2\pi
i/p} .
\]

Then $E_U(f)=\frac1p\sum_{j=0}^{t-1}\sum_{k=0}^{p-1}
(q^{2^j})^k=\frac1p\sum_{j=0}^{t-1} 0 = 0$ since $q^{\alpha}$
is a primitive root of unity when $p$ is prime and $1\leq\alpha<p$.
Likewise, $E_U(f\bar{f}) = \frac1p\sum_{j,j'=0}^{t-1} \sum
_{k=0}^{p-1} q^{k*2^j}\overline{q^{k*2^{j'}}} = \frac1p\sum
_{j=0}^{t-1} p = t$ by the orthogonality relationship of roots of unity.
It follows that $\operatorname{Var}_U(f)=t$.

The Block walk on $\mathbb{Z}_p$ with $n=rt$ steps will be considered where
$r\in\mathbb{N}$ will be chosen later.
Let $P_n(\cdot)$ denote the distribution of
$Z_n=2^{n-1}b_1+2^{n-2}b_2+\cdots+b_n$ induced by $n$ steps of the
block walk starting at state $u=0$. A generic increment will be denoted
by $b$, since the $b_i$ are independent random variables from the same
distribution.

It will be useful to introduce a bit of notation.

If $\alpha\in\mathbb{Z}$ and $y\in\mathbb{Z}_p$, then define
\[
\mu_{\alpha}(y)=\mathsf{\operatorname{Pr}}[2^{\alpha}b=y]=\mathsf
{\operatorname{Pr}}[b=y2^{-\alpha}] .
\]

Recall the Fourier transform of a distribution $\nu$ on $\mathbb
{Z}_p$ is
given by
$\hat{\nu}(\ell)=\sum_{k=0}^{p-1}q^{\ell k}\nu(k)=E_{\nu}
(q^{\ell k})$.
Properties of certain Fourier transforms are required in our work.
First, since $p=2^t-1$ then $\mu_\alpha(y)=\mu_{\alpha+ct}(y)$ for
$c\in\mathbb{N}$,
and so $\hat\mu_\alpha(y)=\hat\mu_{(\alpha \operatorname{mod} t)}(y)$.
By this and (\ref{eqn:convolution}), if $y\in\mathbb{Z}_p$, then
\[
\hat{P_n}(y)=\prod_{\alpha=0}^{n-1}\hat\mu_{\alpha}(y)=
\Biggl(\prod_{\alpha=0}^{t-1}\hat\mu_{\alpha}(y)\Biggr)^r=\hat
{P_t}(y)^r .
\]
Also, $\mu_\alpha(2^jy) = \mu_{(\alpha-j)}(y)$, and so
\[
\forall j\in\mathbb{N}\dvtx \hat{P_t}(2^j y)=\prod_{\alpha=0}^{t-1}
\hat
\mu_{\alpha}(2^j y) = \prod_{\alpha=0}^{t-1} \hat\mu_{\alpha
}(y) = \hat{P_t}(y) .
\]
Finally, for $0\leq j\leq t-1$, define
\[
\Pi_j = \hat{P_t}(2^j-1) = \prod_{\alpha=0}^{t-1} \hat\mu_{\alpha
}(2^j-1)  .
\]
Note that $\hat{P_t}(1-2^j)=\hat{P_t}(2^j-1)$ because $\mu_\alpha
(y)=\mu_\alpha(-y)$ when the step sizes are $\{1,x\}=\{1,-1\}$. Also,
$\Pi_{t-j}=\Pi_j$ because, modulo $p$, $2^{t-j}-1=2^{-j}-1$, and so
$\hat{P_t}(2^{t-j}-1)=\hat{P_t}(2^{-j}-1)=\hat
{P_t}(2^j(2^{-j}-1))=\hat{P_t}(2^j-1)$.

Now turn to mean and variance:
\begin{eqnarray*}
E_{P_n}(f) &=& \sum_k P_n(k) \sum_{j=0}^{t-1} q^{k2^j}
= \sum_{j=0}^{t-1} \hat{P_n}(2^j)
 \\
&=& \sum_{j=0}^{t-1} \hat{P_t}(2^j)^r
= \sum_{j=0}^{t-1} \hat{P_t}(1)^r
= t \Pi_1^r.
\end{eqnarray*}
Likewise,
\begin{eqnarray*}
E_{P_n}(f\bar{f}) &=& \sum_k P_n(k) \sum_{j,j'=0}^{t-1}q^{k(2^j-2^{j'})}
= \sum_{j,j'=0}^{t-1} \hat{P_n}(2^j-2^{j'})
\\
&=& \sum_{j,j'=0}^{t-1} \hat{P_t}(2^j-2^{j'})^r
= \sum_{j=0}^{t-1} \sum_{\beta=0}^{t-1} \hat{P_t}\bigl(2^j(1-2^\beta)\bigr)^r
\\
&=& t \sum_{\beta=0}^{t-1} \hat{P_t}(2^\beta-1)^r
= t \sum_{\beta=0}^{t-1}\Pi_\beta^r.
\end{eqnarray*}

It follows that
\[
\operatorname{Var}_{P_n}(f)=E_{P_n}(f\bar{f})-E_{P_n}(f)E_{P_n}(\bar{f})
= t \sum_{j=0}^{t-1} \Pi_j^r - t^2|\Pi_1|^{2r}.
\]

To apply these relations in Chebyshev's inequality the quantities $\Pi
_j$ need to be examined further.
Let $b$ denote the increment taken but with arithmetic NOT done modulo
$p$, so that $b=p+1$ is possible, that is, repeatedly decide with
probability $1/3$ whether to terminate and if not then flip a coin and
decide whether to add $+1$ or $-1$.
Let $a_k=\mathsf{\operatorname{Pr}}[b=k]$, and note that also $a_k=\mathsf
{\operatorname{Pr}}[b=-k]$ since the
nondoubling steps are symmetric, that is, $u\to u+1$ and $u\to u+x=u-1$.
Then $a_k$ satisfies the recurrence relation
\[
a_k = \tfrac13 (a_{k-1}+a_{k+1}), \qquad a_0=\tfrac13+\tfrac23
a_1,\qquad a_{\infty}=0
\]
which has solution $a_k=\frac{1}{\sqrt{5}}(\frac{3-\sqrt
{5}}{2})^{|k|}$.
For $y\in\mathbb{R}$, let
\[
G(y)
= \sum_{k=-\infty}^{\infty} a_k e^{2\pi iky}
= \frac{1}{\sqrt{5}}
\frac{1-({(3-\sqrt{5})/2})^2}
{1+({(3-\sqrt{5})/2})^2-(3-\sqrt{5})\cos(2\pi y)}.
\]
Since
\[
\mu_\alpha(k)=\mathsf{\operatorname{Pr}}[b=k2^{-\alpha}]=\mathsf{\operatorname{Pr}}[b2^{\alpha
}\equiv k
\operatorname{mod} p]=\sum_{\{b: b2^\alpha\equiv k \operatorname{mod} p\}} a_b,
\]
then
\begin{eqnarray*}
\Pi_j
&=& \prod_{\alpha=0}^{t-1} \hat\mu_\alpha(2^j-1)
= \prod_{\alpha=0}^{t-1} \sum_{k=-\infty}^{\infty} q^{k(2^j-1)}\mu
_\alpha(k)
 \\
&=& \prod_{\alpha=0}^{t-1} \sum_{b=-\infty}^{\infty} a_b e^{2\pi
i b2^{\alpha}(2^j-1)/p}
\\
&=& \prod_{\alpha=0}^{t-1} G\biggl(\frac{2^{\alpha
}(2^j-1)}{p}\biggr) .
\end{eqnarray*}

We can now show the necessary bounds. Recall that $n=rt$ for some
sequence $r=r(t)\in\mathbb{N}$ to be defined.
Let $\lambda=\lambda(t)\in\mathbb{R}$ be a sequence such that
$\lambda
\stackrel{t\to\infty}{\hbox to 20pt{\rightarrowfill}}\infty$, $\lambda=o(\log t)$, and
\[
r = \frac{\log t}{2\log(1/|\Pi_1|)}-\lambda
\]
is an integer. Such a sequence will exist if $|\Pi_1|$ is bounded away
from $0$.

\begin{Claim}\label{c1} $|\Pi_1|$ is bounded away from $0$ and $1$ as
$p=2^t-1\to\infty $.
\end{Claim}

\begin{pf}
First, a few preliminary calculations are necessary. If $y\in[0,1/2]$, then
\[
G(y)\geq
\frac{1}{\sqrt{5}} \frac{1-({(3-\sqrt{5})/2}
)^2}{1+({(3-\sqrt{5})/2})^2+(3-\sqrt{5})}
= \frac15 .
\]
Also, since $y\in\mathbb{R}$, then $|G(y)|\leq\sum a_k|e^{2\pi i k
y}| =
\sum a_k=1$ with equality at $y=0$.
And, since the first derivative satisfies the relation $0 \geq G'(y)
\geq-2\pi(3-\sqrt{5}) > -5$,
then $G(\varepsilon) \geq G(0)-5\varepsilon=1-5\varepsilon$ for $\varepsilon
\geq0$.

To show that $\Pi_1$ is bounded away from $1$ observe that
\[
\Pi_1 \leq G(2^{t-1}/p)*1^{t-1} \stackrel{t\to\infty}{\hbox to 20pt{\rightarrowfill}
} G(1/2) =
\tfrac15.
\]

To bound $\Pi_1$ away from $0$ note that since $G(y)$ is decreasing
for $y\in[0,1/2]$, then
\begin{eqnarray*}
\Pi_1
&=& \prod_{\alpha=0}^{t-1} G(2^{\alpha}/p)
\\
&=& \prod_{\alpha=t-5}^{t-1} G(2^{\alpha}/p) \prod_{\alpha
=0}^{t-6} G(2^{\alpha}/p)
\\
&\geq& \prod_{\beta=1}^5 \frac15 \prod_{\beta=6}^{t}
(1-5\cdot2^{-\beta}-5\cdot2^{-t})
\\
&\geq& 5^{-5}e^{-{5/(2^5(1-5/2^6-5/2^t))}}.
\end{eqnarray*}
The final inequality used the relation is $\ln(1-x) \geq\frac{-x}{1-x}$.
\end{pf}

\begin{Claim}
$(\operatorname{Var}_U(f))^{1/2}=o\bigl(|E_{P_n}(f)-E_U(f)|\bigr)$.
\end{Claim}

\begin{pf}
Since $|\Pi_1|$ is bounded away from $0$ and $1$, then
\[
|E_{P_n}(f) - E_U(f)| = t |\Pi_1|^r = \sqrt{t} |\Pi_1|^{-\lambda}.
\]
The claim follows as $(\operatorname{Var}_U(f))^{1/2}=\sqrt{t}$ and
$|\Pi
_1|^\lambda\to0$.
\end{pf}

\begin{Claim} $(\operatorname{Var}_{P_n}(f))^{1/2} =
o\bigl(|E_{P_n}(f)-E_U(f)|\bigr)$.
\end{Claim}

\begin{pf}
Assume
\[
\frac1t \sum_{j=0}^{t-1}\biggl(\frac{\Pi_j}{\Pi_1^2}
\biggr)^r\stackrel{t\to\infty}{\hbox to 20pt{\rightarrowfill}} 1 .
\]
Then the claim follows from
\[
(\operatorname{Var}_{P_n}(f))^{1/2} = t|\Pi_1|^r\sqrt
{\frac1t \sum
_{j=0}^{t-1} \biggl(\frac{\Pi_j}{\Pi_1^2}\biggr)^r - 1}
= o(t|\Pi_1|^r) .
\]

Hildebrand \cite{Hil061} requires $4$ pages (pages 351--354) to prove
the assumption, albeit for a different function $G(y)$.
Fortunately we do not need to rework his proof as it does not make
explicit use of $G(y)$ but instead depends on only a few properties of it.
There are two facts required in his proof. First, he shows the following:

\begin{Fact} There is some $t_0$ such that if $t\geq t_0$,
then $\Pi_j\leq\Pi_1$ for all $j\geq1$.
\end{Fact}

Hildebrand's proof of this utilizes the following properties:
$G(y)=G(1-y)=G(-y)$, $G(y)$ is decreasing when $y\in[0,1/4]$, and
$G(\frac12-2^{-i})$ is decreasing in \mbox{$i\geq3$} [this corrects an error
where he should have claimed
\[\lim_{p=2^t-1\to\infty} G\biggl(\frac
{2^{t-1}-2^{t-j-1}}{p}\biggr)\leq G(3/8)
\]
 instead of $G(1/2)$]. These
conditions apply to our choice of $G(y)$ as is shown in the proof of
our Claim \ref{c1}.

The second necessary tool is the following:

\begin{Fact} There exists constants $c_0,t_1$ such that
for $t\geq t_1$ and $t^{1/3}\leq j\leq t/2$, then
\[
\frac{\Pi_j}{\Pi_1^2}\leq1+\frac{c_0}{2^j}.
\]
\end{Fact}

Hildebrand's proof utilizes the following properties: $G'(0)=0$,
$|G'(y)|\leq A$, and $G(y)\geq B$ for all $y$ and some $A,B>0$.
All three conditions were shown in the proof of our Claim \ref{c1}, and so his
argument will carry through.

The two facts can now be combined to finish the proof. Let $t$ be
sufficiently large that $\frac{c_0}{2^{(t^{1/3})}}\leq1/r$ which is
possible as $r=O(\log t)$. Recalling that $\Pi_{t-j}=\Pi_j$, then
\begin{eqnarray*}
\sum_{j=1}^{t-1} \biggl|\biggl(\frac{\Pi_j}{\Pi_1^2}
\biggr)^r-1\biggr|
&\leq& 2\sum_{1\leq j<t^{1/3}} \frac{1}{\Pi_1^r}
+ 2\sum_{t^{1/3} \leq j\leq t/2} r^2\frac{c_0}{2^j}
 \\
&\leq& 2t^{1/3} t^{1/2}\Pi_1^{\lambda} + 2\frac t2\frac
{r^2c_0}{2^{(t^{1/3})}}
\\
&=& o(t).
\end{eqnarray*}
Since $\Pi_0=1$, then
\[
\frac1t\sum_{j=0}^{t-1} \biggl(\frac{\Pi_j}{\Pi_1^2}\biggr)^r
\leq\frac1t\biggl(\frac{1}{\Pi_1^{2r}} + (t-1) + o(t)\biggr) = 1 +
o(1) .
\]

In the other direction, since $\operatorname{Var}_{P_n}(f)\geq0$,
then $\frac
1t\sum_{j=0}^{t-1} (\frac{\Pi_j}{\Pi_1^2})^r\geq1$.
\end{pf}
\noqed\end{pf*}
\end{appendix}

\section*{Acknowledgments}
The authors thank S.~Kijima, S.~Miller, I.~Mironov, R.~Venkatesan and
D.~Wilson for several helpful
discussions.

%********************** References ********************************
% imsref loaded by elazauskaite, 2009-10-08 11:06:53
%

\printaddresses


\begin{thebibliography}{22}

%%b1 ###
%%
%(\byear{1986}).
%%

%b2 ###
\bibitem{AF}
%
\begin{bmisc}[auto:springertagbib-v1.0]
\bauthor{\bsnm{Aldous},~\bfnm{David}\binits{D.}} \AND
\bauthor{\bsnm{Fill},~\bfnm{J.}\binits{J.}}
(\byear{2009}).
\textit{Reversible Markov Chains and
Random Walks on Graphs}.  Available at
\href{http://www.stat.berkeley.edu/users/aldous/RWG/book.html}{http://www.stat.berkeley.edu/users/aldous/RWG/book.}
\href{http://www.stat.berkeley.edu/users/aldous/RWG/book.html}{html}.
\end{bmisc}
%
\endbibitem

%b3 ###
\bibitem{CDG871}
%
\begin{barticle}[mr]
\bauthor{\bsnm{Chung},~\bfnm{F.~R.~K.}\binits{F.~R.~K.}},
\bauthor{\bsnm{Diaconis},~\bfnm{Persi}\binits{P.}} \AND
\bauthor{\bsnm{Graham},~\bfnm{R.~L.}\binits{R.~L.}}
(\byear{1987}).
\btitle{Random walks arising in random number generation}.
\bjournal{Ann. Probab.}
\bvolume{15}
\bpages{1148--1165}.
\bid{mr={893921}}
\end{barticle}
%
\endbibitem

%b4 ###
\bibitem{CP05}
%
\begin{bbook}[mr]
\bauthor{\bsnm{Crandall},~\bfnm{Richard}\binits{R.}} \AND
\bauthor{\bsnm{Pomerance},~\bfnm{Carl}\binits{C.}}
(\byear{2005}).
\btitle{Prime Numbers}, \bedition{2nd} ed.
\bpublisher{Springer}, \baddress{New York}.
\bid{mr={2156291}}
\end{bbook}
%
\endbibitem

%%b5 ###
%%
%(\byear{1991}).
%{M}arkov chains, with an application to the exclusion process}.
%%

%b6 ###
\bibitem{Hil061}
%
\begin{barticle}[mr]
\bauthor{\bsnm{Hildebrand},~\bfnm{Martin}\binits{M.}}
(\byear{2006}).
\btitle{On the {C}hung--{D}iaconis--{G}raham random process}.
\bjournal{Electron. Comm. Probab.}
\bvolume{11}
\bpages{347--356 (electronic)}.
\bid{mr={2274529}}
\end{barticle}
%
\endbibitem

%b7 ###
\bibitem{GR91}
%
\begin{barticle}[mr]
\bauthor{\bsnm{Le~Gall},~\bfnm{Jean-Fran{\c{c}}ois}\binits{J.-F.}}
\AND
\bauthor{\bsnm{Rosen},~\bfnm{Jay}\binits{J.}}
(\byear{1991}).
\btitle{The range of stable random walks}.
\bjournal{Ann. Probab.}
\bvolume{19}
\bpages{650--705}.
\bid{mr={1106281}}
\end{barticle}
%
\endbibitem

%b8 ###
\bibitem{KMT071}
%
\begin{bincollection}[auto:springertagbib-v1.0]
\bauthor{\bsnm{Kim},~\bfnm{J-H.}\binits{J-H.}}
\bauthor{\bsnm{Montenegro},~\bfnm{R.}\binits{R.}}
\AND
\bauthor{\bsnm{Tetali},~\bfnm{P.}\binits{P.}}
(\byear{2007}).
\btitle{Near
optimal bounds for collision in Pollard Rho for discrete log}. In
\bbooktitle{Proceedings of the 48th Annual IEEE Symposium on Foundations of Computer
Science} \bpages{215--223}.
\bpublisher{IEEE Computer Society Press},
\baddress{Los Alamitos, CA}.%
\end{bincollection}
%
\endbibitem

%b9 ###
\bibitem{LPS03}
%
\begin{barticle}[mr]
\bauthor{\bsnm{Lyons},~\bfnm{Russell}\binits{R.}},
\bauthor{\bsnm{Peres},~\bfnm{Yuval}\binits{Y.}} \AND
\bauthor{\bsnm{Schramm},~\bfnm{Oded}\binits{O.}}
(\byear{2003}).
\btitle{Markov chain intersections and the loop-erased walk}.
\bjournal{Ann. Inst. H. Poincar\'e Probab. Statist.}
\bvolume{39}
\bpages{779--791}.
\bid{mr={1997212}}
\end{barticle}
%
\endbibitem

%%b10 ###
%%
%(\byear{1989}).
%chains---A combinatorial treatment of expanders}.
%In \bbooktitle{Proceedings of the 30th
%Annual IEEE Symposium on Foundations of Computer Science}
%%

%b11 ###
\bibitem{MV061}
%
\begin{bincollection}[mr]
\bauthor{\bsnm{Miller},~\bfnm{Stephen~D.}\binits{S.~D.}} \AND
\bauthor{\bsnm{Venkatesan},~\bfnm{Ramarathnam}\binits{R.}}
(\byear{2006}).
\btitle{Spectral analysis of {P}ollard Rho collisions}.
In \bbooktitle{Algorithmic Number Theory}.
\bseries{Lecture Notes in Computer Science}
\bvolume{4076}
\bpages{573--581}.
\bpublisher{Springer}, \baddress{Berlin}.
\bid{mr={2282950}}
\end{bincollection}
%
\endbibitem

%b12 ###
\bibitem{MV081}
%
\begin{barticle}[mr]
\bauthor{\bsnm{Miller},~\bfnm{Stephen~D.}\binits{S.~D.}} \AND
\bauthor{\bsnm{Venkatesan},~\bfnm{Ramarathnam}\binits{R.}}
(\byear{2009}).
\btitle{Non-degeneracy of {P}ollard Rho collisions}.
\bjournal{Int. Math. Res. Not. IMRN}
\bvolume{1}
\bpages{1--10}.
\bid{mr={2471294}}
\end{barticle}
%
\endbibitem

%b13 ###
\bibitem{NP08}
%
\begin{bmisc}[auto:springertagbib-v1.0]
\bauthor{\bsnm{Nazarov},~\bfnm{F.}\binits{F.}}
\AND
\bauthor{\bsnm{Peres},~\bfnm{Y.}\binits{Y.}}
 (2010). The birthday problem for
finite reversible Markov chains: A~uniform bound. Preprint.
\end{bmisc}
%
\endbibitem

%b14 ###
\bibitem{Pak021}
%
\begin{bincollection}[auto:springertagbib-v1.0]
\bauthor{\bsnm{Pak},~\bfnm{I.}\binits{I.}}
 (\byear{2002}).
 \btitle{Mixing time and long paths in graphs}. In
\bbooktitle{Proceedings of the 13th Annual ACM-SIAM Symposium on Discrete Algorithms}
\bpages{321--328}.
\bpublisher{ACM},
\baddress{New York}.
\end{bincollection}
%
\endbibitem

%b15 ###
\bibitem{PohHel78}
%
\begin{barticle}[mr]
\bauthor{\bsnm{Pohlig},~\bfnm{Stephen~C.}\binits{S.~C.}} \AND
\bauthor{\bsnm{Hellman},~\bfnm{Martin~E.}\binits{M.~E.}}
(\byear{1978}).
\btitle{An improved algorithm for computing logarithms over {${\rm
GF}(p)$} and
its cryptographic significance}.
\bjournal{IEEE Trans. Information Theory}
\bvolume{IT-24}
\bpages{106--110}.
\bid{mr={0484737}}
\end{barticle}
%
\endbibitem

%b16 ###
\bibitem{Pol75}
%
\begin{barticle}[mr]
\bauthor{\bsnm{Pollard},~\bfnm{J.~M.}\binits{J.~M.}}
(\byear{1975}).
\btitle{A {M}onte {C}arlo method for factorization}.
\bjournal{Nordisk Tidskr. Informationsbehandling (BIT)}
\bvolume{15}
\bpages{331--334}.
\bid{mr={0392798}}
\end{barticle}
%
\endbibitem

%b17 ###
\bibitem{Pol78}
%
\begin{barticle}[mr]
\bauthor{\bsnm{Pollard},~\bfnm{J.~M.}\binits{J.~M.}}
(\byear{1978}).
\btitle{Monte {C}arlo methods for index computation {$({\rm mod} p)$}}.
\bjournal{Math. Comp.}
\bvolume{32}
\bpages{918--924}.
\bid{mr={0491431}}
\end{barticle}
%
\endbibitem

%b18 ###
\bibitem{Pom}
%
\begin{bincollection}[mr]
\bauthor{\bsnm{Pomerance},~\bfnm{Carl}\binits{C.}}
(\byear{2008}).
\btitle{Elementary thoughts on discrete logarithms}.
In \bbooktitle{Algorithmic Number Theory: Lattices, Number Fields,
Curves and
Cryptography}.
\bseries{Mathematical Sciences Research Institute Publications}
\bvolume{44}
\bpages{385--396}.
\bpublisher{Cambridge Univ. Press}, \baddress{Cambridge}.
\bid{mr={2467551}}
\end{bincollection}
%
\endbibitem

%b19 ###
\bibitem{Shoup97}
%
\begin{bincollection}[mr]
\bauthor{\bsnm{Shoup},~\bfnm{Victor}\binits{V.}}
(\byear{1997}).
\btitle{Lower bounds for discrete logarithms and related problems}.
In \bbooktitle{Advances in Cryptology---{EUROCRYPT}'97 ({K}onstanz)}.
\bseries{Lecture Notes in Computer Science}
\bvolume{1233}
\bpages{256--266}.
\bpublisher{Springer}, \baddress{Berlin}.
\bid{mr={1603068}}
\end{bincollection}
%
\endbibitem

%b20 ###
\bibitem{Tes98}
%
\begin{bincollection}[mr]
\bauthor{\bsnm{Teske},~\bfnm{Edlyn}\binits{E.}}
(\byear{1998}).
\btitle{Speeding up {P}ollard's Rho method for computing discrete logarithms}.
In \bbooktitle{Algorithmic Number Theory ({P}ortland, {OR}, 1998)}.
\bseries{Lecture Notes in Computer Science}
\bvolume{1423}
\bpages{541--554}.
\bpublisher{Springer}, \baddress{Berlin}.
\bid{mr={1726100}}
\end{bincollection}
%
\endbibitem

%b21 ###
\bibitem{Tes01}
%
\begin{bincollection}[mr]
\bauthor{\bsnm{Teske},~\bfnm{Edlyn}\binits{E.}}
(\byear{2001}).
\btitle{Square-root algorithms for the discrete logarithm problem (a survey)}.
In \bbooktitle{Public-key Cryptography and Computational Number Theory
({W}arsaw, 2000)}
\bpages{283--301}.
\bpublisher{de Gruyter}, \baddress{Berlin}.
\bid{mr={1881641}}
\end{bincollection}
%
\endbibitem

%b22 ###
\bibitem{VW99}
%
\begin{barticle}[mr]
\bauthor{\bparticle{van }\bsnm{Oorschot},~\bfnm{Paul~C.}\binits
{P.~C.}} \AND
\bauthor{\bsnm{Wiener},~\bfnm{Michael~J.}\binits{M.~J.}}
(\byear{1999}).
\btitle{Parallel collision search with cryptanalytic applications}.
\bjournal{J. Cryptology}
\bvolume{12}
\bpages{1--28}.
\bid{mr={1664774}}
\end{barticle}
%
\endbibitem

\end{thebibliography}
\end{document}